\documentclass[12pt,a4paper]{article}
\usepackage{bbm}
\usepackage{mathrsfs}
\usepackage{graphicx}
\usepackage{amsmath}
\usepackage{amssymb}
\usepackage{amsfonts}
\usepackage{enumerate}
\usepackage{theorem}

\makeatletter
\def\tank#1{\protected@xdef\@thanks{\@thanks
        \protect\footnotetext[0]{#1}}}
\def\bigfoot{

    \@footnotetext}
\makeatother

\newcommand{\ea}{\end{array}}
\newtheorem{theorem}{Theorem}[section]

\newtheorem{lemma}{Lemma}[section]
\newtheorem{definition}{Definition}[section]
\newtheorem{Rem}{Remark}[section]

{\theorembodyfont{\rmfamily}
}

\title
{\bf The asymptotic behavior of primitive equations with multiplicative noise \thanks{This work was partially
supported by NNSF of China(Grant No. 11401057),   Natural Science Foundation Project of CQ  (Grant No. cstc2016jcyjA0326),
Fundamental Research Funds for the Central Universities(Grant No. 106112015CDJXY100005) and China Scholarship Council (Grant No.:201506055003).} }
\author{
Rangrang Zhang
\thanks{
Beijing Institute of Technology, No. 5 Zhongguancun South
Street, Haidian District, Beijing, 100081, P R China }
\tank{rrzhang@amss.ac.cn.}
\qquad
 Guoli Zhou
\thanks{ Chongqing University,  No.174 Shazheng Street, Shapingba District, Chongqing, 400044, P R China }
\tank{E-mail:zhouguoli736@126.com.}}

\begin{document}
\maketitle

\begin{abstract}
This paper is concerned with the existence of invariant measure for 3D stochastic primitive equations driven by linear multiplicative noise under non-periodic boundary conditions. The common method is to apply Sobolev imbedding theorem  to proving the tightness of the distribution of the solution. However, this method fails because of the non-linearity and non-periodic boundary conditions of the stochastic primitive equations. To overcome the difficulties, we show the existence of random attractor by proving the compact property and the regularity of the solution operator. Then we show the existence of invariant measure.
\end{abstract}



\section{Introduction}
\par
The paper is concerned with the stochastic primitive equations  (PEs) in a cylindrical domain $$\mathcal{O}=M\times (-h,0)\subset \mathbb{R}^{3},$$ where $M$ is  a smooth bounded domain in $\mathbb{R}^{2}$:
\begin{eqnarray}\label{z-1}
\partial_{t} \upsilon+(\upsilon\cdot \nabla)\upsilon+ w \partial_{z} \upsilon +f \upsilon^{\perp} +\nabla p+L_{1}\upsilon&=&\sum_{k=1}^{n}\alpha_{k}\upsilon \circ d w_{k}^{1},\\
\label{z-2}
\partial_{z}p +T&=&0,\\
\label{z-3}
\nabla\cdot \mathbf{\upsilon}+\partial_{z} w&=&0,\\
\label{z-4}
\partial_{t} T+\upsilon\cdot \nabla T+w \partial_{z} T+L_{2}T&=&Q+ \sum_{k=1}^{n}\beta_{k}T \circ d w_{k}^{2}.
\end{eqnarray}
The unknowns for the 3D stochastic PEs are the fluid velocity field $(\upsilon,w )=(\upsilon_{1},\upsilon_{2},w )\in \mathbb{R}^{3}$ with  $\upsilon=(\upsilon_{1},\upsilon_{2})$ and $ \upsilon^{\perp}=(-\upsilon_{2}, \upsilon_{1} ) $ being horizontal, the temperature $T$ and the pressure $p.$
$f=f_{0}(\beta+y)$ is the given Coriolis parameter, $Q$ is a given heat source.
$\nabla = (\partial_{x}, \partial_{y})$ is the horizontal gradient operator and $\Delta=\partial_{x}^{2}+ \partial_{y}^{2}$
is the horizontal Laplacian. Let $\alpha_{i},\beta_{i} \in \mathbb{R}, i=1,2,\cdot\cdot\cdot,n$, $\{w_{i}^{1},w_{i}^{2},i=1,2...n\}$ be
a sequence of one-dimensional, independent, identically distributed  Brownian
motions.   Here, we take $\sum_{k=1}^{n}\alpha_{k}\upsilon \circ dw_{k}^{1}$ and $\sum_{k=1}^{n}\beta_{k}T \circ d w_{k}^{2}$ to be  Stratonovich multiplicative noise.

 The viscosity and the heat diffusion operators $L_{1}$ and $L_{2}$ are given by
\begin{eqnarray*}
L_{i}=-\nu_{i}\Delta-\mu_{i}\partial_{zz} , \quad  i=1,2,
\end{eqnarray*}
where $\nu_{1}, \mu_{1}$ are the horizontal and vertical Reynolds numbers and
$\nu_{2}, \mu_{2}$ are the horizontal and vertical heat diffusivity. Without loss of generality, we assume that
$$
\nu_{1}= \mu_{1}=\nu_{2}= \mu_{2}=1.
$$
\par
The boundary of $\mathcal{O}$ is partitioned into three parts: $\Gamma_{u}\cup\Gamma_{b}\cup\Gamma_{s}, $ where
\begin{eqnarray*}
&&\Gamma_{u}=\{(x,y,z)\in \overline{\mathcal{O}}:z=0 \},\\
&&\Gamma_{b}=\{(x,y,z)\in \overline{\mathcal{O}}:z=-h \},\\
&&\Gamma_{s}=\{(x,y,z)\in \overline{\mathcal{O}}:(x,y)\in \partial M, -h \leq z\leq 0 \}.
\end{eqnarray*}
Here $h$ is a sufficiently smooth function. Without loss of generality, we assume $h=1.$

We impose the following boundary conditions on the 3D stochastic PEs.
\begin{eqnarray}\label{z-5}
\partial_{z} \upsilon =\eta,\ \ w=0,\ \ \partial_{z} T=-\gamma (T-\tau)\ \ &{\rm{on}}& \  \Gamma_{u},\\
\label{z-6}
\partial_{z} \upsilon=0,\ \ w=0,\ \ \partial_{z} T=0\ \ &{\rm{on}}& \ \Gamma_{b},\\
\label{z-7}
\upsilon\cdot \vec{n}=0,\ \ \partial_{\vec{n}}\upsilon\times \vec{n} =0,\ \  \partial_{\vec{n}} T=0\ \ &{\rm{on}}&\  \Gamma_{s},
\end{eqnarray}
where $\eta(x,y)$ is the wind stress on the surface of the ocean, $\gamma$ is a positive constant, $\tau$ is the typical temperature distribution on the top surface of the ocean and $\vec{n}$ is the norm vector to $\Gamma_{s}.$ For the sake of simplicity, we assume that $Q$ is independent of time and $\eta=\tau=0$. It is worth pointing out that all results obtained in this paper are still valid for the general case by making some simple modifications.

Integrating (\ref{z-3}) from $-1$ to $z$ and using  (\ref{z-5}), (\ref{z-6}), we have
\begin{equation}
w(t,x,y,z):=\Phi(\upsilon )(t,x,y,z)=-\int^{z}_{-1}\nabla\cdot \upsilon (t,x,y,z')dz',
\end{equation}
moreover,
\[
\int^{0}_{-1}\nabla\cdot \upsilon  dz=0.
\]
Integrating (\ref{z-2}) from $-1$ to $z$, set $p_{s}$ be a certain unknown function at $\Gamma_{b}$ satisfying
\[
p(x,y,z,t)= p_{s}(x,y,t)-\int^{z}_{-1} T(x,y,z',t) dz'.
\]
Then (\ref{z-1})-(\ref{z-4}) can be rewritten as
\begin{eqnarray}\label{e-1}
&\partial_{t} \upsilon+L_{1}\upsilon+(\upsilon\cdot \nabla)\upsilon+\Phi(v) \partial_{z}\upsilon
+\nabla p_{s}-\int_{-1}^{z}\nabla Tdz' +f \upsilon^{\perp} =\sum_{k=1}^{n}\alpha_{k}\upsilon \circ d w_{k}^{1},&\\
\label{e-2}
&\partial_{t} T+L_{2}T+\upsilon\cdot \nabla T+\Phi(v) \partial_{z} T=Q+\sum_{k=1}^{n}\beta_{k}T \circ d w_{k}^{2},&\\
\label{e-3-1}
&\int^{0}_{-1}\nabla\cdot \upsilon  dz=0.&
\end{eqnarray}
The boundary value and initial conditions for (\ref{e-1})-(\ref{e-3-1}) are given by
\begin{eqnarray}\label{e-3}
\partial_{z} \upsilon|_{\Gamma_{u}}=\partial_{z} \upsilon|_{\Gamma_{b}}=0,\quad
\upsilon\cdot \vec{n}|_{\Gamma_{s}}=0, \quad \partial_{\vec{n}}\upsilon\times \vec{n}|_{\Gamma_{s}}=0,\\
\label{e-4}
\Big{(}\partial_{z}T+\gamma T\Big{)}|_{\Gamma_{u}}=\partial_{z}T|_{\Gamma_{b}}=0, \quad \partial_{\vec{n}}T|_{\Gamma_{s}}=0;\\
\label{e-5}
\upsilon(x,y,z,t_0)=\upsilon_{0}(x,y,z), \ \ T(x,y,z,t_0)=T_{0}(x,y,z).
\end{eqnarray}

\par
The primitive equations are the basic model used in the study of climate and
weather prediction, which can be used to describe the motion of the atmosphere when the hydrostatic
assumption is enforced (see \cite{G,H1,H2} and the references therein). This model has been intensively investigated because of the interests stemmed from physics and mathematics. As far as we know, their mathematical study was initiated by Lions, Temam and Wang (see e.g. \cite{LTW1,LTW2,LTW3,LTW4}). For example, the existence of global weak solutions for the primitive equations was established in \cite{LTW2}. Guill\'{e}n-Gonz\'{a}lez et al. obtained the global existence of strong solutions to the primitive equations with small initial data in \cite{GMR}. The local existence of strong solutions to the primitive equations under the small depth hypothesis was established by Hu et al. in \cite{HTZ}. Cao and Titi developed a beautiful approach to dealing with the $L^6$-norm of the fluctuation $\tilde{v}$ of horizontal velocity and obtained the global well-posedness for the 3D viscous primitive equations in \cite{CT1}. Subsequently, in \cite{KZ}, Kukavica and Ziane developed a different method to handle non-rectangular domains and boundary conditions with physical reality. The existence of
the global attractor was established by Ju in \cite{J}. For the global well-posedness of 3D primitive equations with partial dissipation, we refer the reader to some papers, see e.g.  $\cite{CIN, CLT1, CLT2, CLT3, CT2}$.

\par
Along with the great successful developments of deterministic
primitive equations, the random case has also been developed rapidly. In \cite{GH}, Guo and Huang obtained some kind of weak type compactness properties of the strong solution under the condition that the momentum equation is driven by an additive stochastic forcing and the
thermodynamical equation is driven by a deterministic heat source. Debussche, Glatt-Holtz, Temam and Ziane concerned the global well-posedness of strong solution when the primitive
equations are driven by multiplicative stochastic forcing in \cite{D-G-T-Z}.  The ergodic theory and large deviations for the 3D stochastic primitive equations were studied by Dong, Zhai and Zhang in \cite{DZZ1,DZZ2}.
When the primitive equations are driven by an infinite-dimensional fractional noise taking values in some Hilbert space, the latter author of this paper obtained the existence of random attractor in \cite{Z}.
Under the periodic conditions, Glatt-Holtz, Kukavica, Vicol and Ziane constructed an invariant measure for the 3D PEs in \cite{GKVZ},
 where some delicate and careful  techniques are involved. However, it seems that in their paper, the conditions on the noise don't cover the linear multiplicative noise.

\par
In this article, we aim to prove the existence of invariant measures for 3D stochastic PEs driven by linear multiplicative noise under non-periodic boundary conditions. The common method is to apply Krylov-Bogoliubov lemma (see \cite{DZ}) which requires a uniform estimates with respect to the initial data in $H^2(\mathcal{O})$. However, it is quite difficult because of highly non-linear drift terms and non-periodic boundary conditions. Instead of using that method, we provide a new method to find a compact absorbing set in $H^{1}(\mathcal{O})$, which guarantees the existence of random attractor in $H^{1}(\mathcal{O})$. The main idea is that we firstly prove the solution operator of 3D stochastic PEs is compact in $H^{1}(\mathcal{O})$, $\mathbb{P}-$a.s.. Then we can construct a compact absorbing ball in $H^{1}(\mathcal{O})$ by the compact property of solution operator. It's important to point out that the random attractor we obtained is stronger than that in \cite{GH}. Specifically, the random attractor we obtain is compact in $H^{1}(\mathcal{O})$, $\mathbb{P}$-a.s. and attracts any orbit starting from $-\infty$ in the strong topology of $H^{1}(\mathcal{O})$ while the random attractor in \cite{GH} is not necessarily a compact subset in  $H^{1}(\mathcal{O})$. In addition, the linear multiplicative noise as well as  non-periodic boundary conditions are handled at the same time in this paper, which are different from  \cite{GKVZ}.

%
Taking into account the asymptotical compact property of the solution operator, we can prove the existence of invariant measures by showing that the one-point motions associated with the flow generated by 3D PEs define a family of Markov processes. Up to now, there are no works on the existence of invariant measures for the stochastic PEs subject to non-periodic boundary conditions. Maybe it is an attempt to solve this problem by proving the asymptotic compact property of the solution operator.

\par
The remaining of this paper is organized as follows. In section 2,  some preliminaries of 3D stochastic primitive equations are stated.  In section 3, the global well-posedness of 3D stochastic primitive equations is proved. In section 4, we establish the existence of random attractor. Finally, in section 5, the existence of invariant measures for 3D stochastic primitive equations is obtained. As usual, constants
$C$ may change from one line to the next, unless, we give a special declaration. Denote by $C(a)$ a constant which depends on some parameter $a$.
\section{Preliminaries}
For $1\leq p\leq \infty,$ let $L^{p}(\mathcal{O}), L^{p}(M)$ be the usual Lebesgue spaces with the norm
\begin{eqnarray*}
 |\phi|_p=\left\{
            \begin{array}{ll}
              \left(\int_{\mathcal{O}}|\phi(x,y,z)|^pdxdydz\right)^{\frac{1}{p}}, & \phi\in L^p(\mathcal{O}), \\
               \left(\int_M|\phi(x,y)|^pdxdy\right)^{\frac{1}{p}}, &  \phi\in L^p(M).
            \end{array}
          \right.
\end{eqnarray*}
In particular, $|\cdot|$ and $(\cdot,\cdot)$ represent norm and inner product of $L^2(\mathcal{O})$ (or $L^2(M)$), respectively. For $m\in \mathbb{N}_+$, $(W^{m,p}(\mathcal{O}), \|\cdot\|_{m,p})$ stands for the classical Sobolev space, see \cite{ARA}. When $p=2$, we denote by $H^m(\mathcal{O})=W^{m,2}(\mathcal{O})$ with norm $\|\cdot\|_{m}$.
Without confusion, we shall sometimes abuse notation and denote by $\|\cdot\|_{m}$ the norm in $ H^{m}(M)$.  Let
\begin{eqnarray*}
\mathcal{V}_{1}&=&\Big\{\mathbf{\upsilon}\in (C^{\infty}(\mathcal{O}))^{2}: \partial_{z} {\upsilon} |_{z=0}=0, \ \partial_{z}\upsilon|_{z=-1}=0,\ \mathbf{\upsilon}\cdot \vec{n}|_{\Gamma_{s}}=0,\\
 &&\ \partial_{\vec{n}} \upsilon \times \vec{n}|_{\Gamma_{s}}=0, \ \int_{-1}^{0}\nabla \cdot \mathbf{\upsilon} dz=0 \Big\},\\
\mathcal{V}_{2}&=&\Big\{T\in C^{\infty}(\mathcal{O}): \partial_{z} T |_{z=-1}=0, \ (\partial_{z} T+\gamma T)|_{z=0}=0,\  \partial_{\vec{n}} T |_{\Gamma_{s}}=0\Big \}.
\end{eqnarray*}
We denote by $V_{1}$ the closure space of $\mathcal{V}_{1}$ in $(H^{1}(\mathcal{O}))^{2}$ and $V_{2}$ the closure space of $\mathcal{V}_{2}$ in $H^{1}(\mathcal{O})$, respectively. Let $H_{1}$ be the closure space of $\mathcal{V}_{1}$ with respect to the norm $|\cdot|_{2}.$ Define $H_{2}=L^{2}(\mathcal{O}).$ Set
$$V= V_{1}\times V_{2}, \quad H=H_{1}\times H_{2} .$$
Let $U=(\mathbf{\upsilon},T)$, $\tilde{U}=(\tilde{\mathbf{\upsilon}},\tilde{T})$, $V$ is equipped with the inner product
\begin{eqnarray*}
\langle U,  \tilde{U}\rangle_{V}&:=&\langle \mathbf{\upsilon},  \tilde{\mathbf{\upsilon}}\rangle_{V_{1}}+\langle T,  \tilde{T}\rangle_{V_{2}},\\
\langle \mathbf{\upsilon},  \tilde{\mathbf{\upsilon}}\rangle_{V_{1}}&:=&\int_{\mathcal{O}}\left(\nabla v\cdot \nabla \tilde{v}+\frac{\partial v}{\partial z}\cdot\frac{\partial \tilde{v}}{\partial z}\right)dxdydz,\\
\langle T,  \tilde{T}\rangle_{V_{2}}&:=&\int_{\mathcal{O}}\left(\nabla T\cdot \nabla \tilde{T}+\frac{\partial T}{\partial z}\frac{\partial \tilde{T}}{\partial z}\right)dxdydz+\gamma\int_{\Gamma_{u}} T_{1}T_{2}dxdy.
\end{eqnarray*}
Subsequently, the norm of $V$ is defined by $\|U\|=\langle U,  U\rangle_{V}^{\frac{1}{2}}.$ The inner product of $H$ is defined by
\begin{eqnarray*}
\langle U,  \tilde{U}\rangle_{{H}}&:=&\langle \mathbf{\upsilon},  \tilde{\mathbf{\upsilon}}\rangle+\langle T,  \tilde{T}\rangle,\\
\langle \mathbf{\upsilon},  \tilde{\mathbf{\upsilon}}\rangle&:=&\int_{\mathcal{O}} \mathbf{\upsilon}\cdot  \tilde{\mathbf{\upsilon}}dxdydz,\\
\langle T,  \tilde{T}\rangle&:=&\int_{\mathcal{O}} T  \tilde{T}dxdydz.
\end{eqnarray*}
Denote $V_{i}'$ the dual space of $V_{i}$, $i=1,2.$ 
Furthermore, we have the compact embedding relationship
\begin{eqnarray*}
D(A_{i})\subset V_{i}\subset H_{i}\subset V_{i}'\subset D(A_{i})',
\end{eqnarray*}
and
\begin{eqnarray*}
\langle \cdot,  \cdot \rangle_{V_{i}}=\langle A_{i}\cdot, \cdot    \rangle= \langle A_{i}^{\frac{1}{2}}\cdot, A_{i}^{\frac{1}{2}}\cdot\rangle, \ \ i=1,2.
\end{eqnarray*}
For the sake of simplicity, in the following, we denote
\[
\int_{\mathcal{O}}\cdot \ dxdydz=\int_{\mathcal{O}}\cdot \ , \quad \int_{M}\cdot\  dxdy=\int_{M}\cdot \ .
\]

\section{Global well-posedness of (\ref{e-1})-(\ref{e-5})}
In this section, we aim to prove the global well-posedness of (\ref{e-1})-(\ref{e-5}). Firstly, we introduce the following definition.
Given $\mathcal{T}>0,$ fix a single stochastic basis $(\Omega, \mathcal{F}, \{\mathcal{F}_{t_{0},t}\}_{t\in[t_0,\mathcal{T}]}, \mathbb{P})$, where
\begin{eqnarray}\label{e-13}
\mathcal{F}_{t_{0},t}:=\sigma(W_{k}^{j}(s)-W_{k}^{j}(t_{0}) , s\in [t_{0}, t], j=1,2).
\end{eqnarray}
\begin{definition}
Fix $\mathcal{T}>0$, a continuous $V-$valued $\mathcal{F}_{t_{0},t}-$adapted random field $(U(.,t))_{t\in [t_0,\mathcal{T}]}=(\upsilon(.,t),T(.,t))_{t\in [t_{0},\mathcal{T}]}$ defined on $(\Omega, \mathcal{F}, \mathbb{P})$ is said to be a strong (weak) solution  to (\ref{e-1})-(\ref{e-5}) if
\[
U\in C([t_0,\mathcal{T}];V)\cap L^{2}([t_0,\mathcal{T}];(H^{2}(\mathcal{O}))^{3})\ ( U\in C([t_0,\mathcal{T}];{H})\cap L^{2}([t_0,\mathcal{T}];(H^{1}(\mathcal{O}))^{3})\quad \mathbb{P}-a.s..
\]
and the following
\begin{eqnarray*}
&&\int_{\mathcal{O}}\!\upsilon(t)\!\cdot  \phi_{1} -\int^{t}_{t_0}ds\int_{\mathcal{O}}\{[(\upsilon\cdot \nabla )\phi_{1} +\Phi(\upsilon)\partial_{z}\phi_{1} ]\upsilon-[(fk\times \upsilon)\cdot \phi_{1}+(\int_{-1}^{z}Tdz')\nabla \cdot \phi_{1}] \}\nonumber\\
&&+ \ \int_{t_0}^{t}ds\int_{\mathcal{O}}\upsilon\cdot L_{1}\phi_{1} =\int_{\mathcal{O}}\upsilon_{0}\cdot \phi_{1}+\int_{t_{0}}^{t}\int_{\mathcal{O}}\sum_{k=1}^{n}\alpha_{k} \upsilon \circ d w^{1}_{k}(s,w)\cdot \phi_{1},\nonumber\\
&&\int_{\mathcal{O}}T(t)\phi_{2}-\int_{t_0}^{t}ds\int_{\mathcal{O}}[(\upsilon\cdot \nabla)\phi_{2}+ \Phi(\upsilon)\partial_{z}\phi_{2}]T +\int_{t_0}^{t}ds\int_{\mathcal{O}}T L_{2}\phi_{2}\\
&=&\int_{\mathcal{O}}T_{0}\phi_{2}+\int_{t_0}^{t}ds\int_{\mathcal{O}}Q\phi_{2} +\int_{t_{0}}^{t}\int_{\mathcal{O}}\sum_{k=1}^{n}\beta_{k} T \circ d w^{2}_{k}(s,w)\cdot \phi_{2},
\end{eqnarray*}
holds $\mathbb{P}-$a.s., for all $t\in [t_0,\mathcal{T}]$ and $\phi=(\phi_{1},\phi_{2})\in D(A_{1})\times D(A_{2}).$
\end{definition}
Consider
\begin{eqnarray*}
\alpha(t)=\exp(- \sum_{k=1}^{n}\alpha_{k} w_{k}^{1}),\ \beta(t)=\exp(- \sum_{k=1}^{n}\beta_{k} w_{k}^{2}).
\end{eqnarray*}
Then $\alpha (t)$ and $\beta(t)$ satisfy the following stratonovich equations
\begin{eqnarray*}
d\alpha(t)=-\sum_{k=1}^{n}\alpha_{k}\alpha(t)\circ dw_{k}^{1}(t),\ \
d\beta(t)=-\sum_{k=1}^{n}\beta_{k}\beta(t)\circ dw_{k}^{2}(t).
\end{eqnarray*}
Define
\begin{eqnarray*}
(u(t), \theta(t))=(\alpha(t)v(t),\ \beta(t)T(t)).
\end{eqnarray*}
Then, $(u(t), \theta(t))$ satisfies
\begin{eqnarray}\label{qq-1}
&\partial_{t} u-\Delta u-\partial_{zz} u+\alpha^{-1} u\cdot\nabla u+\alpha^{-1}\Phi(u) \partial_{z}u+f u^{\bot} +\alpha\nabla p_{s}-\alpha\beta^{-1}\int_{-1}^{z}\nabla \theta dz'=0,&\\
\label{qq-2}
&\partial_{t}\theta -\Delta \theta-\partial_{zz}\theta + \alpha^{-1}u\cdot\nabla \theta+\alpha^{-1}\Phi(u) \partial_{z} \theta=\beta Q,&\\
\label{qq-3}
&\int^{0}_{-1}\nabla\cdot udz=0.&
\end{eqnarray}
The boundary and initial conditions for (\ref{qq-1})-(\ref{qq-3}) are
\begin{eqnarray}\label{qq-4}
\partial_{z} u|_{\Gamma_{u}}=\partial_{z}u|_{\Gamma_{b}}=0,
u\cdot \vec{n}|_{\Gamma_{s}}=0, \partial_{\vec{n}}u\times \vec{n}|_{\Gamma_{s}}=0,\\
\label{qq-5}
\Big{(}\partial_{z}\theta+\gamma \theta\Big{)}|_{\Gamma_{u}}=\partial_{z}\theta|_{\Gamma_{b}}=0, \ \ \partial_{\vec{n}}\theta|_{\Gamma_{s}}=0,\\
\label{qq-6}
(u\big{|}_{t_0}, \theta\big{|}_{t_0})=(\upsilon_{0}, T_{0}).
\end{eqnarray}
\begin{definition}\label{d-1}
Let $\mathcal{T}$ be a fixed positive time and $(v_{0}, T_{0})\in V$.  $ (u,\theta)$ is called a strong solution of the system (\ref{qq-1})-(\ref{qq-6}) on the time interval $[t_0,\mathcal{T}]$ if it satisfies (\ref{qq-1})-(\ref{qq-6}) in the weak sense such that $\mathbb{P}$-a.s.
\begin{eqnarray*}
&&u\in C([t_0,\mathcal{T}];V_{1})\cap L^{2}([t_0,\mathcal{T}]; (H^{2}(\mathcal{O}))^{2}),\\
&&\theta\in C([t_0,\mathcal{T}];V_{2})\cap L^{2}([t_0,\mathcal{T}]; H^{2}(\mathcal{O})).
\end{eqnarray*}
\end{definition}

\begin{theorem}
[\textbf{Existence of local solutions to (\ref{qq-1})-(\ref{qq-6})}] If $Q\in L^{2}(\mathcal{O}), v_{0}\in V_{1}$, $T_{0}\in V_{2}$. Then, for $\mathbb{P}$-a.s., $\omega \in\Omega,$ there exists a stopping time $T^{*}>0$
such that $(u,\theta)$ is a strong solution of the system (\ref{qq-1})-(\ref{qq-6}) on the interval $[t_0, T^{*}].$
\end{theorem}

%

The proof of  the existence of local solutions to (\ref{qq-1})-(\ref{qq-6}) is similar to \cite{GMR}, we omit it. Before showing the global well-posedness of the strong solution, we recall the following Lemma, a special case of a general result of Lions and Magenes $\cite{LM}$,  which will help us to show the continuity of the solution with respect to time in $(H^{1}(\mathcal{O}))^{3}.$ For the proof, see $\cite{T2}.$
\begin{lemma}\label{l-1}
Let $V , H, V'$ be three Hilbert spaces such that $V \subset H = H' \subset V' $
, where $ H'$ and $V'$ are the dual spaces of $H$ and $V$, respectively. Suppose $u \in
L^{2}(0, T; V )$ and $\frac{du}{dt}\in L^{2}(0, T; V')$. Then $u$ is almost everywhere equal to a continuous function from $[0, T]$ into $H$.
\end{lemma}

\begin{theorem}\label{thm-2}
{\rm{[\textbf{Existence of global solution to (\ref{qq-1})-(\ref{qq-6})}]}} If $Q\in L^{2}(\mathcal{O}), v_{0}\in V_{1}$, $T_{0}\in V_{2}$, and $\mathcal{T}>0$. Then, for $\mathbb{P}$-a.s., $\omega \in\Omega,$ there exists a unique strong solution $(u, \theta)$ of the system (\ref{qq-1})-(\ref{qq-6}) or equivalently $(v,T)$ to the system (\ref{e-1})-(\ref{e-5}) on the interval $[t_0, \mathcal{T}]$.
\end{theorem}

\begin{flushleft}
\textbf{Proof.} \quad Let $[t_0,\tau_{*})$ be the maximal interval of existence of the strong solution, for fixed $\omega\in \Omega,$
we will establish various norms of this solution in the interval $[t_0, \tau_{*}).$
In particular, we will show that if $\tau_{*}<\infty$ then $H^1-$norm of the strong solution is
bounded over the interval $[t_0, \tau_{*})$.

\textbf{A priori estimates:} \quad
 Referring to \cite{CT1}, define
\begin{eqnarray*}
\bar{\phi}(x,y)=\int_{-1}^{0}\phi(x,y,\xi)d\xi,\ \ \forall\ (x,y)\in M.
\end{eqnarray*}
In particular,
\begin{eqnarray*}
\bar{u}(x,y)=\int_{-1}^{0}u(x,y,\xi)d\xi,\ \ \ \mathrm{in}\ M.
\end{eqnarray*}
Let
\begin{eqnarray*}
\tilde{u}=u-\bar{u}.
\end{eqnarray*}
Notice that
\begin{eqnarray*}
\bar{\tilde{u}}=0, \quad \nabla \cdot \bar{u}=0 \ \ \mathrm{in}\ M.
\end{eqnarray*}
Taking the average of equations (\ref{qq-1}) in the $z$ direction
over the interval $(-1,0)$, and using boundary conditions (\ref{qq-4}), we have
\begin{eqnarray}\label{qq-41}
\partial_{t} \bar{u}+\alpha^{-1}\overline{u\cdot\nabla u+\Phi(u)\partial_{z}u }
+f \bar{u}^{\bot}+\alpha\nabla p_{s}-\alpha \beta^{-1}\int_{-1}^{0}\int_{-1}^{z}\nabla \theta dz'dz-\Delta \bar{u}=0.
\end{eqnarray}
By the integration by parts, we get
\begin{eqnarray}\label{qq-42}
\int_{-1}^{0}\Phi(u)\partial_{z}u dz&=&\int_{-1}^{0}u \nabla \cdot u dz=
\int_{-1}^{0}(\nabla\cdot\tilde{u})
\tilde{u}dz,\\
\label{qq-43}
\int_{-1}^{0}u\cdot \nabla u dz&=& \int_{-1}^{0}\tilde{u}\cdot \nabla \tilde{u}dz
+ \bar{u}\cdot \nabla \bar{u}.
\end{eqnarray}
Substituting (\ref{qq-42}) and (\ref{qq-43}) into (\ref{qq-41}), $\bar{u}$ satisfies
\begin{eqnarray}\notag
& \partial_{t} \bar{u}-\Delta \bar{u}+\alpha^{-1}(\overline{\tilde{u}\cdot \nabla\tilde{u}}
+\overline{\tilde{u} \nabla\cdot\tilde{u}}+\bar{u}\cdot \nabla \bar{u})+f\bar{u}^{\bot}+\alpha\nabla p_{s}&\\
\label{qq-20}
& -\alpha\beta^{-1}\nabla \int_{-1}^{0}
\int_{-1}^{z}\theta(x,y,\lambda,t)d\lambda dz=0,&\\
\label{qq-21}
& \nabla \cdot \bar{u}=0\ \ \ {\rm{in}}\ M,&\\
\label{qq-22}
& \bar{u}\cdot \vec{n} =0,\ \partial_{\vec{n}} \bar{u} \times \vec{n}=0\  {\rm{on}}\  M.&
\end{eqnarray}
By subtracting (\ref{qq-20}) from (\ref{qq-1}) and using (\ref{qq-4}), (\ref{qq-22}), we conclude that $\tilde{u}$ satisfies
\begin{eqnarray}\notag
& \partial_{t} \tilde{u}-\Delta \tilde{u}- \partial_{zz}\tilde{u}+\alpha^{-1}\tilde{u}\cdot
 \nabla\tilde{u}+\alpha^{-1}\Phi(\tilde{u}) \partial_{z} \tilde{u}
 +\alpha^{-1}\tilde{u}\cdot \nabla\bar{u}+\alpha^{-1}\bar{u}\cdot \nabla\tilde{u}&
\\
 \label{qq-15}
  &\ - \alpha^{-1}\overline{\tilde{u}\cdot \nabla \tilde{u}}-\alpha^{-1}\overline{\tilde{u} \nabla\cdot\tilde{u}}+f\tilde{u}^{\bot}-\alpha \beta^{-1}\int_{-1}^{z}\nabla \theta dz'+
 \alpha\beta^{-1} \int_{-1}^{0}  \int_{-1}^{z}\nabla \theta dz'dz=0, &\\
 \label{qq-16}
& \partial_{z} \tilde{u}|_{z=0}=0,\ \  \partial_{z} \tilde{u}|_{z=-1}=0,\ \  \tilde{u}\cdot \vec{n}|_{\Gamma_{s}}=0,\ \
\partial_{ \vec{n}} \tilde{u}\times \vec{n}|_{\Gamma_{s}}=0.&
\end{eqnarray}
In the following, we will study the properties of $\bar{u}$ and $\tilde{u}.$

\textbf{(1)\quad Estimates of $|\theta|^2$ and $|u|^2$.}\quad
Take the inner product of equation (\ref{qq-2}) with $\theta$ in $H_{2}$, we get
\begin{eqnarray*}
&&\frac{1}{2}\partial_{t} |\theta|^{2}+|\nabla \theta|^{2}+|\theta_{z}|^{2}+\gamma|\theta(z=0)|^{2}\\
&=&\beta\int_{\mathcal{O}}Q\theta -\alpha^{-1}\int_{\mathcal{O}}\Big{(}u\cdot\nabla \theta+\Phi (u) \partial_{z} \theta\Big{)}\theta .
\end{eqnarray*}
By integration by parts,
\begin{eqnarray*}
\alpha^{-1}\int_{\mathcal{O}}\Big{(}u\cdot\nabla \theta+\Phi (u)\partial_{z} \theta \Big{)}\theta =0.
\end{eqnarray*}
By the H\"{o}lder inequality, we have
\begin{eqnarray*}
\frac{1}{2}\partial_{t} |\theta|^{2}+|\nabla \theta|^{2}+|\theta_{z}|^{2}+\gamma|\theta(z=0)|^{2}
\leq \beta\int_{\mathcal{O}}Q\theta \leq \varepsilon |\theta|^2+C\beta^2|Q|^2,
\end{eqnarray*}
referring to (48) in \cite{CT1},
\begin{eqnarray*}
|\theta|^{2}\leq 2 |\partial_{z} \theta|^{2}+2|\theta(z=0)|^{2},
\end{eqnarray*}
then, we arrive at
\begin{eqnarray}\label{qq-8}
\partial_{t} |\theta|^{2}+2|\nabla \theta|^{2}+(2-4\varepsilon)|\theta_{z}|^{2}+(2\gamma-4\varepsilon)|\theta(z=0)|^{2}
\leq  C\beta^{2}|Q|^{2}.
\end{eqnarray}
Hence, there exists a positive $\lambda$ such that
\begin{eqnarray*}
\partial_{t} |\theta|^{2}+\lambda| \theta|^{2}
\leq  C\beta^{2}|Q|^{2}.
\end{eqnarray*}
Applying the Gronwall inequality, we have
\begin{eqnarray}\label{qq-9}
|\theta(t)|^{2}
\leq  |\theta_{t_{0}}|^{2} e^{-\lambda(t-t_{0} )}+C\int_{t_{0}}^{t}\beta^{2}e^{\lambda(s-t)}|Q|^{2}ds.
\end{eqnarray}
In view of (\ref{qq-8}) and (\ref{qq-9}), we obtain
\begin{eqnarray}\label{qq-10}
\sup\limits_{t\in [t_{0}, \tau^{*})}|\theta(t)|^{2}+ \int_{t_{0}}^{\tau^{*}}\|\theta(t)\|^{2}dt
\leq  C.
\end{eqnarray}
Taking inner product of (\ref{qq-1}) with $u$ in $H_{1}$, by integration by parts, we have
\begin{eqnarray}\label{r-10}
\partial_{t} |u|^{2}+(1-\varepsilon)|\nabla u|^{2}+|u_{z}|^{2}
\leq  C|\theta|^{2}.
\end{eqnarray}
By (\ref{qq-10}), we get
\begin{eqnarray}\label{qq-26}
\sup\limits_{t\in [t_{0}, \tau^{*})}|u(t)|^{2}+ \int_{t_{0}}^{\tau^{*}}\|u(t)\|^{2}dt
\leq  C.
\end{eqnarray}
\textbf{(2)\quad Estimates of $|\theta|^4_{4}$ and  $|\tilde{u}|^4_{4}$.} \quad
Taking the inner product of the equation (\ref{qq-2}) with $\theta^{3}$ in $H_{2}$, we have
\begin{eqnarray}
&&\frac{1}{4}\partial_{t}|\theta|^{4}_{4}+\frac{3}{4}|\nabla \theta^{2}|^{2}+\frac{3}{4}|(\theta^{2})_{z}|^{2}+\gamma\int_{M}|\theta(z=0)|^{4}\nonumber\\
&=&\beta\int_{\mathcal{O}}Q\theta^{3}-\alpha^{-1}\int_{\mathcal{O}}[u\cdot \nabla \theta+\Phi( u)\partial_{z} \theta ]\theta^{3}.
\end{eqnarray}
By integration by parts, we have
\begin{eqnarray}\label{qq-13}
\alpha^{-1}\int_{\mathcal{O}}[u\cdot \nabla \theta+\Phi( u)\partial_{z} \theta ]\theta^{3}=0.
\end{eqnarray}
Applying the interpolation inequality to $|\theta^{2}|_{3}$, we obtain
\begin{eqnarray}\label{qq-11}
|\theta^{2}|_{3}\leq C |\theta^{2}|^{\frac{1}{2}}(|\nabla \theta^{2}|^{\frac{1}{2}}+ |\partial_{z} \theta^{2}|^{\frac{1}{2}}+\alpha|\theta^{2}(z=0) |^{\frac{1}{2}}   ).
\end{eqnarray}
Using (\ref{qq-11}) and the H\"{o}lder inequality, we get
\begin{eqnarray}\label{qq-12}
\int_{\mathcal{O}}\beta Q\theta^{3}\leq \varepsilon (|\nabla \theta^{2}|^{2}+   |\partial_{z} \theta^{2}|^{2}+\alpha|\theta^{2}(z=0) |^{2} )+C\beta|Q|^{\frac{8}{5}}|\theta|_{4}^{\frac{12}{5}}.
\end{eqnarray}
Combining (\ref{qq-13}) and (\ref{qq-12}), we arrive at
\begin{eqnarray}\label{qq-14}
\partial_{t}|\theta|^{4}_{4}+|\nabla \theta^{2}|^{2}+|(\theta^{2})_{z}|^{2}+\alpha\int_{M}|\theta(z=0)|^{4}
\leq  C\beta|Q|^{\frac{8}{5}}|\theta|_{4}^{\frac{12}{5}}.
\end{eqnarray}
Since
\[
\theta^4(x,y,z)=-\int^0_z \partial_r \theta^4(x,y,r)dr+\theta^4(z=0),
\]
by the Young's inequality, we have
\begin{eqnarray*}
|\theta|_{4}^{4}&=&-\int_{\mathcal{O}}\int_{z}^{0}\partial_{r}\theta^{4}dr +\int_{M}\int_{-1}^{0}\theta^{4}(z=0)dz\\
&\leq&\frac{1}{2}|\theta|_{4}^{4}+ 8 |\partial_{z}(\theta^{2})|^{2}+\int_{M}\theta^{4}(z=0),
\end{eqnarray*}
then
\begin{eqnarray*}
|\theta|_{4}^{4}\leq 16|\partial_{z}(\theta^{2}) |^{2}+2|\theta(z=0)|_{4}^{4}.
\end{eqnarray*}
From (\ref{qq-14}), we get
\begin{eqnarray*}
\partial_{t}|\theta|^{4}_{4}+|\theta|^{4}_{4}&\leq& C\beta|Q|^{\frac{8}{5}}|\theta|_{4}^{\frac{12}{5}},\\
\partial_{t} |\theta|^{2}_{4}+|\theta|^{2}_{4}
&\leq& C\beta|Q|^{\frac{8}{5}}|\theta|_{4}^{\frac{2}{5}}.
\end{eqnarray*}
Applying the Gronwall inequality, there exists a positive number which is still denoted by $\lambda$ such that
\begin{eqnarray}\label{e-12}
|\theta(t)|^{2}_{4}\leq |\theta_{t_{0}}|^{2}_{4}e^{-\lambda(t-t_{0})}+C\int_{t_0}^{t}\beta(s) e^{-\lambda(t-s)}|Q|^{2}ds
\end{eqnarray}
for $t\in[t_0, \tau_{*}).$

Taking the inner product of the equation (\ref{qq-15}) with $|\tilde{u}|^{2}\tilde{u}$ in $H_{1}$, we obtain
\begin{eqnarray*}
&&\frac{1}{4}\partial_{t} |\tilde{u}|_{4}^{4}+
\frac{1}{2}\int_{\mathcal{O}}\Big{(}|\nabla(|\tilde{u} |^{2}) |^{2}+ |\partial_{z}(|\tilde{u} |^{2}) |^{2}\Big{)}
+\int_{\mathcal{O}}|\tilde{u}|^{2}(|\nabla \tilde{u} |^{2}+|\partial_{z}\tilde{u} |^{2})\nonumber\\
&=& -\alpha^{-1}\int_{\mathcal{O}}((\tilde{u}\cdot \nabla)\tilde{u}+\Phi(\tilde{u})
\partial_{z} \tilde{u} )\cdot |\tilde{u}|^{2}\tilde{u}\\
&&\ -\alpha^{-1}\int_{\mathcal{O}}(\tilde{u}\cdot \nabla \bar{u})\cdot |\tilde{u}|^{2}\tilde{u} -\alpha^{-1}\int_{\mathcal{O}}(\bar{u}\cdot \nabla \tilde{u})\cdot |\tilde{u}|^{2}\tilde{u}\\
&&\ +\alpha^{-1}\int_{\mathcal{O}}\overline{\tilde{u}\nabla \cdot \tilde{u}
+\tilde{u}\cdot \nabla \tilde{u} }\cdot| \tilde{u}|^{2}\tilde{u}\\
&&\ +\alpha\beta^{-1}\int_{\mathcal{O}}(\int^z_{-1}\nabla \theta dz'-\int^0_{-1}\int^z_{-1}\nabla \theta dz'dz)\cdot| \tilde{u}|^{2}\tilde{u}.
\end{eqnarray*}
By integration by parts and boundary conditions (\ref{qq-16}), we have
\begin{eqnarray*}
&&\int_{\mathcal{O}}((\tilde{u}\cdot \nabla)\tilde{u}+\Phi(\tilde{u})
\partial_{z} \tilde{u} )|\tilde{u}|^{2}\tilde{u}=0,
\\
&&\int_{\mathcal{O}}(\bar{u} \cdot\nabla \tilde{u})|\tilde{u}|^{2}\tilde{u}=-\frac{1}{4}\int_{\mathcal{O}}|\tilde{u}|^{4}\nabla
\cdot \bar{u}=0,
\end{eqnarray*}
and
\begin{eqnarray*}
\int_{\mathcal{O}}[(\tilde{u}\cdot \nabla )\bar{u} ]\cdot |\tilde{u}|^{2}\tilde{u}&=&-\int_{\mathcal{O}}[(\tilde{u} \cdot \nabla )
|\tilde{u}|^{2}\tilde{u}]\cdot \bar{u}-\int_{\mathcal{O}} (\nabla \cdot \tilde{u}  ) |\tilde{u}|^{2}\tilde{u} \cdot \bar{u},\\
\int_{\mathcal{O}}\overline{\tilde{u}\nabla \cdot \tilde{u}
+\tilde{u}\cdot \nabla \tilde{u} }\cdot| \tilde{u}|^{2}\tilde{u}&=&-\int_{\mathcal{O}}\overline{\tilde{u}_{k} \tilde{u}_{j} }\partial_{x_{k}}(| \tilde{u}|^{2}\tilde{u}_{j} ),
\end{eqnarray*}
where $\tilde{u}_{k}$ is the $k-$th coordinate of  $\tilde{u}$, $k=1,2.$

Based on the above equalities and  by integration by parts, we obtain
\begin{eqnarray}\notag
&&\frac{1}{4}\partial_{t} |\tilde{u}|_{4}^{4}+
\frac{1}{2}\int_{\mathcal{O}}\Big{(}|\nabla(|\tilde{u} |^{2}) |^{2}+ |\partial_{z}(|\tilde{u} |^{2}) |^{2}\Big{)}
+\int_{\mathcal{O}}|\tilde{u}|^{2}(|\nabla \tilde{u} |^{2}+|\partial_{z}\tilde{u} |^{2})\\ \notag
&=&\alpha^{-1}\int_{\mathcal{O}}\bar{u} \cdot(\tilde{u} \cdot\nabla )|\tilde{u}|^{2}\tilde{u}
+\alpha^{-1}\int_{\mathcal{O}}(\nabla\cdot \tilde{u} )\bar{u} \cdot
|\tilde{u}|^{2}\tilde{u}\\ \notag
&&\ -\alpha^{-1}\int_{\mathcal{O}}\overline{ \tilde{u}_{k}  \tilde{u}_{j}} \partial_{x_{k}}(| \tilde{u}|^{2}\tilde{u}_{j} )-\alpha\beta^{-1}\int_{\mathcal{O}}\Big{(}\int_{-1}^{z}\theta d\lambda- \int_{-1}^{0}\int_{-1}^{z}\theta d\lambda dz  \Big{)}\nabla \cdot |\tilde{u}|^{2}\tilde{u}\\
\label{qq-17}
&:=&I_1+I_2+I_3+I_4.
\end{eqnarray}
Applying the H\"{o}lder inequality, the Minkowski inequality and the interpolation inequalities,  we obtain
\begin{eqnarray*}
I_{1}
&\leq&\alpha^{-1}\int_{M}|\bar{u}|\int_{-1}^{0}|\tilde{u}||\nabla \tilde{u}| |\tilde{u}|^{2}dz\\
&\leq&\alpha^{-1}\int_{M}|\bar{u}|\Big{(}\int_{-1}^{0}|\tilde{u}|^{2}|\nabla \tilde{u}|^{2}dz \Big{)}^{\frac{1}{2}}
\Big{(}\int_{-1}^{0}|\tilde{u}|^{4}dz\Big{)}^{\frac{1}{2}}\\
&\leq&\alpha^{-1}|\bar{u}|_{L^{4}(M)}|\nabla(|\tilde{u}|^{2})|\Big{(}\int_{M}(\int_{-1}^{0}|\tilde{u}|^{4} dz)^{2} \Big{)}^{\frac{1}{4}}\\
&\leq&\alpha^{-1}|\bar{u}|_{L^{4}(M)}|\nabla(|\tilde{u}|^{2})|
\Big{(}\int_{-1}^{0}(\int_{M}|\tilde{u}|^{8} )^{\frac{1}{2}}dz\Big{)}^{\frac{1}{2}}\\
&\leq&\alpha^{-1}|\bar{u}|_{L^{4}(M)}|\nabla(|\tilde{u}|^{2})|\Big{(}\int_{-1}^{0}|(|\tilde{u}|^{2}) |_{L^{2}(M)} \|(|\tilde{u}|^{2}) \|_{H^{1}(M)}
 dz\Big{)}^{\frac{1}{2}}\\
&\leq&\alpha^{-1}|\bar{u}|_{L^{4}(M)}|\nabla(|\tilde{u}|^{2})||(|\tilde{u}|^{2})|^{\frac{1}{2}}
[|(|\tilde{u}|^{2})|^{\frac{1}{2}}+|\nabla (|\tilde{u}|^{2}) |^{\frac{1}{2}}+|\partial_{z} (|\tilde{u}|^{2})|^{\frac{1}{2}}].
\end{eqnarray*}
By the Young's inequality and the interpolation inequalities, we get
\begin{eqnarray}
I_{1}&\leq& \varepsilon (|\nabla (|\tilde{u}|^{2} )|^{2}+ |\partial_{z} (|\tilde{u}|^{2} )|^{2})+C(\alpha^{-2}|\bar{u}|_{L^{4}(M)}^{2}+\alpha^{-4}|\bar{u}|_{L^{4}(M)}^{4} )|\tilde{u}|^{4}_{4} \nonumber\\
&\leq& \varepsilon (|\nabla (|\tilde{u}|^{2} )|^{2}+ |\partial_{z} (|\tilde{u}|^{2} )|^{2})
+C(\alpha^{-2}\|u\|_{1}^{2}+\alpha^{-4}|u|^{2}\|u\|_{1}^{2} )|\tilde{u}|^{4}_{4}.
\end{eqnarray}
By the H\"{o}lder inequality, the Minkowski inequality, the interpolation inequality and the Young's inequality, we get
\begin{eqnarray*}
I_2&=&\alpha^{-1}\int_{M}\bar{u}\cdot\int_{-1}^{0}|\nabla \tilde{u}||\tilde{u}|^{3}dz\\
&\leq& \alpha^{-1}\int_{M}|\bar{u}|\Big(\int_{-1}^{0}|\nabla \tilde{u}|^{2}|\tilde{u}|^{2}dz\Big)^{\frac{1}{2}} \Big( \int_{-1}^{0}|\tilde{u}|^{4}dz  \Big)^{\frac{1}{2}} \\
&\leq&\alpha^{-1}| |\nabla \tilde{u}| |\tilde{u}| ||\bar{u}|_{4}\Big( \int_{M}\Big( \int_{-1}^{0}|\tilde{u}|^{4}dz    \Big)^{2}\Big)^{\frac{1}{4}}\\
&\leq&\alpha^{-1}| |\nabla \tilde{u}| |\tilde{u}| ||\bar{u}|_{4}\Big( \int_{-1}^{0}\Big( \int_{M}|\tilde{u}|^{2\times4}   \Big)^{\frac{1}{2}}dz\Big)^{\frac{1}{2}}\\
&\leq&C\alpha^{-1}| |\nabla \tilde{u}| |\tilde{u}| ||\bar{u}|_{4}\Big( \int_{-1}^{0}||\tilde{u}|^{2}|_{L^{2}(M)}(|\nabla|\tilde{u}|^{2}|_{L^{2}(M)} +
|\partial_{z}|\tilde{u}|^{2}|_{L^{2}(M)} + ||\tilde{u}|^{2}|_{L^{2}(M)})dz\Big)^{\frac{1}{2}}\\
&\leq&\varepsilon( | |\nabla \tilde{u}| |\tilde{u}| |^{2} +|\nabla |\tilde{u}|^{2} |^{2}+|\partial_{z}|\tilde{u}|^{2}|^{2} )+C(\alpha^{-2}|\bar{u}|_{4}^{2}+\alpha^{-4}|\bar{u}|_{4}^{4}  )||\tilde{u}|^{2}|^{2}\\
&\leq&\varepsilon( | |\nabla \tilde{u}| |\tilde{u}| |^{2} +|\nabla |\tilde{u}|^{2} |^{2}+|\partial_{z}|\tilde{u}|^{2}|^{2} )+C(\alpha^{-2}\|u\|_{1}^{2}+\alpha^{-4}|u|^{2}\|u\|_{1}^{2}  )|\tilde{u}|_{4}^{4}.
\end{eqnarray*}
Applying the H\"{o}lder inequality, the Minkowski inequality, the interpolation inequality and the Young's inequality, we have
\begin{eqnarray*}
I_{3}&=&-\alpha^{-1}\int_{\mathcal{O}}\overline{\tilde{u}_{k}
\tilde{u}_{j}  }\partial_{x_{k}}(| \tilde{u}|^{2}\tilde{u}_{j} )\nonumber\\
&\leq&\alpha^{-1}\int_{M}\Big{(}\int_{-1}^{0} |\tilde{u}|^{2}dz\Big{)}\Big{(}  \int_{-1}^{0} | \nabla\tilde{u}| |\tilde{u}|^{2}dz  \Big{)}\nonumber\\
&\leq& \alpha^{-1}\Big{(}\int_{\mathcal{O}}| \nabla\tilde{u}|^{2} |\tilde{u}|^{2} \Big{)}^{\frac{1}{2}}
\Big{(} \int_{M} \Big{(}\int_{-1}^{0}|\tilde{u}|^{2}dz\Big{)}^{3}\Big{)}^{\frac{1}{2}}\\
&\leq& \alpha^{-1}\Big{(}\int_{\mathcal{O}}| \nabla\tilde{u}|^{2} |\tilde{u}|^{2} \Big{)}^{\frac{1}{2}}
\Big{(}\int_{-1}^{0}\Big{(}\int_{M}|\tilde{u}|^{6}\Big{)}^{\frac{1}{3}}dz \Big{)}^{\frac{3}{2}} \nonumber\\
&\leq& C\alpha^{-1}\Big{(}\int_{\mathcal{O}}| \nabla\tilde{u}|^{2} |\tilde{u}|^{2} \Big{)}^{\frac{1}{2}}
\Big{(}\int_{-1}^{0}|\tilde{u}|^{\frac{4}{3}}_{L^{4}(M)}\cdot\|\tilde{u}\|^{\frac{2}{3}}_{H^{1}(M)} dz \Big{)}^{\frac{3}{2}}\nonumber\\
&\leq& C\alpha^{-1}\Big{(}\int_{\mathcal{O}}| \nabla\tilde{u}|^{2} |\tilde{u}|^{2} \Big{)}^{\frac{1}{2}}
\Big{(}\int^0_{-1}|\tilde{u}|^{4}_{L^{4}(M)}dz\Big{)}^{\frac{1}{2}}\int^0_{-1}\|\tilde{u}\|_{H^{1}(M)}dz\\
&\leq&\varepsilon \int_{\mathcal{O}}| \nabla\tilde{u}|^{2} |\tilde{u}|^{2}+C\alpha^{-2}\| u\|^{2}_{1}|\tilde{u}|_{4}^{4}.
\end{eqnarray*}
Analogously, we have
\begin{eqnarray*}
I_{4}
&\leq&\alpha\beta^{-1} \Big{(} \int_{\mathcal{O}}|\nabla\tilde{u}|^{2}|\tilde{u}|^{2}\Big{)}^{\frac{1}{2}}
\Big{(}\int_{\mathcal{O}}|\tilde{u}|^{4}\Big{)}^{\frac{1}{4}}\Big{(}\int_{\mathcal{O}}|\theta|^{4}\Big{)}^{\frac{1}{4}}\nonumber\\
&\leq& \varepsilon \int_{\mathcal{O}}|\nabla\tilde{u}|^{2}|\tilde{u}|^{2}+C\alpha^{2}\beta^{-2}| \tilde{u}|_{4}^{2}|\theta|_{4}^{2}.
\end{eqnarray*}
Collecting all the above inequalities, we have
\begin{eqnarray}\label{r-20}
&&\partial_{t} |\tilde{u}|^{4}_{4}+
\int_{\mathcal{O}}\Big{(}|\nabla(|\tilde{u} |^{2}) |^{2}+ |\partial_{z}(|\tilde{u} |^{2}) |^{2}\Big{)}
+\int_{\mathcal{O}}|\tilde{u}|^{2}(|\nabla \tilde{u} |^{2}+|\partial_{z}\tilde{u} |^{2})\nonumber\\
&\leq&C\alpha^{2}\beta^{-2}|\theta|_{4}^{2}|\tilde{u}|_{4}^{2}+C(\alpha^{-2}\|u\|_{1}^{2}+\alpha^{-4}|u|^{2}\|u\|_{1}^{2} )|\tilde{u}|_{4}^{4},
\end{eqnarray}
and
\begin{eqnarray}\label{r-12}
\partial_{t}|\tilde{u}|^{2}_{4}&\leq&C\alpha^{2}\beta^{-2}|\theta|_{4}^{2}+C(\alpha^{-2}\|u\|_{1}^{2}+\alpha^{-4}|u|^{2}\|u\|_{1}^{2} )|\tilde{u}|_{4}^{2}.
\end{eqnarray}
Applying the Gronwall inequality, we conclude that
\begin{eqnarray}\notag
&&\sup\limits_{t\in [t_0, \tau*)}|\tilde{u}(t)|_{4}^{4}+
\int_{t_0}^{\tau*}\int_{\mathcal{O}}\Big{(}|\nabla(|\tilde{u} |^{2}) |^{2}+ |\partial_{z}(|\tilde{u} |^{2}) |^{2}\Big{)}ds\\
&&\ +\int_{t_0}^{\tau*}\int_{\mathcal{O}}|\tilde{u}|^{2}(|\nabla \tilde{u} |^{2}+|\partial_{z}\tilde{u} |^{2})ds
\label{qq-18}
\leq C.
\end{eqnarray}
\textbf{(3) \quad Estimates of $|\nabla \bar{u}|^2$ and $|u_z|^2$.}\quad
Taking the inner product of equation (\ref{qq-20}) with $-\Delta \bar{u}$ in $L^{2}(M)$,  we arrive at
\begin{eqnarray*}
\frac{1}{2}\partial_{t} |\nabla \bar{u} |^{2}+|\Delta \bar{u}|^{2}
&=&-\alpha^{-1}\int_{M}\bar{u}\cdot \nabla \bar{u}\cdot \Delta\bar{u}-\alpha^{-1}\int_{M}\overline{\tilde{u}\cdot \nabla\tilde{u}}\cdot \Delta\bar{u}\nonumber\\
&&\ -\alpha^{-1}\int_{M}\overline{\tilde{u} \nabla\cdot \tilde{u}}\cdot \Delta\bar{u}-\int_{M}(\bar{u}^\perp+\alpha \nabla p_s)\cdot \Delta\bar{u} \\ \notag
&&\ +\alpha \beta^{-1}\int_M \nabla \int_{-1}^{0}\int_{-1}^{z}\theta(x,y,\lambda, t)d\lambda dz\cdot\Delta \bar{u}.
\end{eqnarray*}
By integration by parts and (\ref{qq-21})-(\ref{qq-22})(for more detail, see \cite{CT1} ), we have
\begin{eqnarray*}
\int_{M}\bar{u}^\perp\cdot \Delta \bar{u}=0,\ \ \ \
\int_{M}\nabla p_{s}\cdot\Delta \bar{u}=0,\\
\int_{M}\nabla \int_{-1}^{0}\int_{-1}^{z}\theta(x,y,\lambda, t)d\lambda dz\cdot\Delta \bar{u}=0.
\end{eqnarray*}
Applying the H\"{o}lder inequality and the interpolation inequalities, we obtain
\begin{eqnarray*}
\alpha^{-1}\int_{M}\bar{u}\cdot \nabla \bar{u} \cdot \Delta\bar{u}
&\leq&C\alpha^{-1}|\bar{u}|^{\frac{1}{2}}|\nabla \bar{u}||\Delta \bar{u}|^{\frac{3}{2}}
\nonumber\\
&\leq & \varepsilon |\Delta\bar{u}|^{2}+C\alpha^{-4}|\bar{u}|^{2}|\nabla \bar{u}|^{4}.
\end{eqnarray*}
Using the H\"{o}lder inequality, the Minkowski inequality and the Sobolev imbedding theorem, we have
\begin{eqnarray*}
&&\alpha^{-1}\int_{M}\overline{\tilde{u} \nabla\cdot \tilde{u}}\cdot \Delta\bar{u}
+\alpha^{-1}\int_{M}\overline{\tilde{u} \cdot \nabla \tilde{u}}\cdot \Delta\bar{u}\nonumber\\
&\leq& \alpha^{-1}|\Delta \bar{u}|\Big{(}\int_{\mathcal{O}}|\tilde{u} |^{2}  |\nabla\tilde{u} |^{2}   \Big{)}^{\frac{1}{2}}\nonumber\\
&\leq& \varepsilon |\Delta \bar{u}|^{2}+C\alpha^{-2}\int_{\mathcal{O}}|\tilde{u} |^{2}  |\nabla\tilde{u} |^{2}.
\end{eqnarray*}
From the above inequalities, we conclude that
\begin{eqnarray}\label{qq-25}
\partial_{t} |\nabla \bar{u} |^{2} +|\Delta \bar{u}|^{2}&\leq&C\alpha^{-4}|u|^{2}\|u\|_{1}^{2}|\nabla \bar{u}|^{2}+C\alpha^{-2}\int_{\mathcal{O}}|\tilde{u} |^{2}  |\nabla\tilde{u} |^{2}.
\end{eqnarray}
Therefore, applying Gronwall inequality and (\ref{qq-26}), (\ref{qq-18}), (\ref{qq-25}), we arrive at
\begin{eqnarray}\label{qq-28}
\sup\limits_{t\in[t_0, \tau*)}|\nabla \bar{u} (t)|^{2}\leq
C.
\end{eqnarray}
Denote $u_z=\frac{\partial u}{\partial z}$, from (\ref{qq-1}), we get
\begin{eqnarray}\notag
&\partial_{t} u_{z}-\Delta u_{z}-\partial_{zz} u_{z}+\alpha^{-1}
u\cdot\nabla u_{z}+\alpha^{-1}u_{z}\cdot\nabla u-\alpha^{-1} u_{z}\nabla \cdot u+\alpha^{-1}\Phi (u) u_{zz} &\\
\label{qq-29}
 &\ +fk\times u_{z} -\alpha \beta^{-1} \nabla \theta=0.&
\end{eqnarray}
Taking the inner product of the equation (\ref{qq-29}) with $u_{z}$ in $H_{1}$, we obtain
\begin{eqnarray*}
&&\frac{1}{2}\partial_{t} |u_{z}|^{2}+|\nabla u_{z} |^{2}+| u_{zz} |^{2}\\
&=& -\alpha^{-1}\int_{\mathcal{O}}\Big{(}(u\cdot \nabla )u_{z}+\Phi(u)  u_{zz} \Big{)}\cdot u_{z}-\alpha^{-1}\int_{\mathcal{O}}[(u_{z} \cdot \nabla)u]\cdot u_{z}\\
&&\ +\alpha^{-1}\int_{\mathcal{O}}(\nabla \cdot u) u_{z}\cdot u_{z}-\int_{\mathcal{O}} u^\perp_{z}\cdot u_{z} +\alpha\beta^{-1}\int_{\mathcal{O}}\nabla \theta \cdot u_{z}.
\end{eqnarray*}
By integration by parts, we obtain
\begin{eqnarray*}
\alpha^{-1}\int_{\mathcal{O}}\Big{(}(u\cdot \nabla )u_{z}+\Phi(u)  u_{zz} \Big{)}\cdot u_{z}=0 .
\end{eqnarray*}
Thanks to the H\"{o}lder inequality, the interpolation inequality and the Sobolev imbedding theorem, we reach
\begin{eqnarray*}\notag
 \alpha^{-1}\int_{\mathcal{O}}[(u_{z} \cdot \nabla)u]\cdot u_{z}
&\leq & C\alpha^{-1}\int_{\mathcal{O}} |u||u_{z}||\nabla u_{z}|
\nonumber\\
&\leq& C\alpha^{-1}|\nabla u_{z} ||u|_{4}|u_{z}|_{4}
\nonumber\\
&\leq& C\alpha^{-1}|\nabla u_{z} ||u|_{4}|u_{z}|^{\frac{1}{4}}(|\nabla u_{z} |^{\frac{3}{4}}+ |\partial_{z} u_{z} |^{\frac{3}{4}}+
|u_{z}|^{\frac{3}{4}})\nonumber\\
&\leq&\varepsilon (|\nabla u_{z} |^{2}+ |\partial_{z} u_{z} |^{2})
+C(\alpha^{-8}|\nabla\bar{u}|^{8}+1) |u_{z}|^{2}.
\end{eqnarray*}
Similar to the above, we get
\begin{eqnarray*}
\int_{\mathcal{O}}(\nabla \cdot u) u_{z}\cdot u_{z}
\leq \varepsilon (|\nabla u_{z} |^{2}+ |u_{zz} |^{2})
+C(\alpha^{-8}|\nabla\bar{u}|^{8}+1) |u_{z}|^{2}.
\end{eqnarray*}
Collecting the above inequalities, we have
\begin{eqnarray}\label{qq-31}
\partial_{t} |u_{z}|^{2}+|\nabla u_{z} |^{2}+| u_{zz} |^{2}
\leq
C(\alpha^{-8}|\bar{u}|_{4}^{8}+\alpha^{-8}|\nabla\bar{u}|^{8}+1 )|u_{z}|^{2}.
\end{eqnarray}
Applying Gronwall inequality to (\ref{qq-31}), and by (\ref{qq-28}), we reach
\begin{eqnarray}\label{qq-32}
\sup\limits_{t\in[t_0,\tau*)}|u_{z}(t)|^{2}+\int_{t_0}^{\tau*}(|\nabla u_{z}(s) |^{2}+|u_{zz} (s)|^{2})ds
\leq C.
\end{eqnarray}
\textbf{(4) \quad Estimates of $|\nabla u|^2$ and $|\nabla \theta|^2$.}\quad
Taking the inner product of equation (\ref{qq-1}) with $-\Delta u$ in $H_{1}$, we reach
\begin{eqnarray*}
&&\frac{1}{2}\partial_{t} |\nabla u|^{2}+|\Delta u|^{2}+|\nabla \partial_z u|^2\nonumber\\
&=& \alpha^{-1}\int_{\mathcal{O}}[(u\cdot \nabla ) u+\Phi(u)\partial_z u ]\cdot\Delta u
+\int_{\mathcal{O}}fk\times u\cdot\Delta u\\
&&\ +\alpha \int_{\mathcal{O}}\nabla p_s\cdot\Delta u
-\alpha\beta^{-1}\int_{\mathcal{O}}(\int^z_{-1}\nabla \theta dz')\cdot\Delta u.
\end{eqnarray*}
By the H\"{o}lder inequality, the interpolation inequality and the Sobolev inequality, we have
\begin{eqnarray}
&&\alpha^{-1}\int_{\mathcal{O}}[(u\cdot \nabla ) u ]\cdot\Delta u\nonumber\\
&\leq& C \alpha^{-1}|\Delta u||\nabla u|_{4}|u|_{4}\nonumber\\
&\leq& \alpha^{-1}|\Delta u||\nabla u|^{\frac{1}{4}}(|\Delta u|^{\frac{3}{4}}+ |\nabla u_{z}|_{4}^{\frac{3}{4}}+|\nabla u|^{\frac{3}{4}}  )
|u|_{4}\nonumber\\
&\leq& \varepsilon (|\Delta u|^{2}+|\nabla u_{z}|^{2})
+C(\alpha^{-2}|u|_{4}^{2} +\alpha^{-8}|u|_{4}^{8} )
|\nabla u|^{2}\nonumber\\
&\leq& \varepsilon (|\Delta u|^{2}+|\nabla u_{z}|^{2})+C(\alpha^{-2}|\tilde{u}|_{4}^{2}+\alpha^{-2}|\nabla \bar{u}|^{2} +\alpha^{-8}|\tilde{u}|_{4}^{8}+ \alpha^{-8}|\nabla \bar{u}|^{8})
|\nabla u|^{2}.
\end{eqnarray}
Due to the H\"{o}lder inequality, the Minkowsky inequality, the interpolation inequality and the Sobolev embedding theorem, we get
\begin{eqnarray}\notag
&&\alpha^{-1}\int_{\mathcal{O}} \Phi(u) u_{z}\cdot \Delta u\\ \notag
&\leq&\alpha^{-1}\int_{M}\Big{(}\int_{-1}^{0}|\nabla \cdot u |dz\int_{-1}^{0}|u_{z}|\cdot|\Delta u|dz\Big{)}
\\ \notag
&\leq&\alpha^{-1}|\Delta u|\Big{(}\int_{M}(\int_{-1}^{0}|\nabla \cdot u |   dz )^{4} \Big{)}^{\frac{1}{4}}
\Big{(}\int_{M} (\int_{-1}^{0}| u_{z}|^{2} dz  )^{2} \Big{)}^{\frac{1}{4}}\\ \notag
&\leq& C\alpha^{-1}|\Delta u|\Big{(}\int_{-1}^{0}| \nabla u|_{L^{2}(M)}^{\frac{1}{2}}(| \Delta u|_{L^{2}(M)}^{\frac{1}{2}}
+ | \nabla u|_{L^{2}(M)}^{\frac{1}{2}})  dz    \Big{)}\\ \notag
&&\ \times \Big{(}\int_{-1}^{0} |u_{z}|_{L^{2}(M)}(|\nabla u_{z}|_{L^{2}(M)}+ |u_{z}|_{L^{2}(M)}) dz  \Big{)}^{\frac{1}{2}}\\
\notag
&\leq & \varepsilon (|\Delta u|^{2}+|\nabla u_{z}|^{2}   )
+  C|\nabla u|^{2}(\alpha^{-2}|u_{z}|^{2}+\alpha^{-2}|u_{z}||\nabla u_{z}|\\
\label{r-40}
&&\ +\alpha^{-4}|u_{z}|^{4}+\alpha^{-4}|u_{z}|^{2}|\nabla u_{z}|^{2}).
\end{eqnarray}
We also have
\begin{eqnarray*}
\int_{\mathcal{O}}(fk\times u)\cdot\Delta u=0,\quad \int_{\mathcal{O}}\nabla p_{s}\cdot\Delta u=0.
\end{eqnarray*}
Collecting all the above inequalities, we get
\begin{eqnarray}\notag
&&\partial_{t} |\nabla u|^{2}+|\Delta u|^{2}+|\nabla \partial_z u|^2\\ \notag
&\leq&C(\alpha^{2}\beta^{-2}\|\theta\|_{1}^{2}+ \alpha^{-2} |\tilde{u}|_{4}^{2}+\alpha^{-8}|\tilde{u}|_{4}^{8}+ \alpha^{-2}|\nabla \bar{u}|^{2}+ \alpha^{-8}|\nabla \bar{u}|^{8} \\
\label{r-29}
&&+ \alpha^{-2}|u_{z}|^{2}+ \alpha^{-4} |u_{z}|^{4}+ \alpha^{-2}|\nabla u_{z}|^{2}+\alpha^{-4} |u_{z}|^{2}|\nabla u_{z}|^{2})
|\nabla u|^{2}.
\end{eqnarray}
Applying the Gronwall inequality, and by (\ref{qq-10}),(\ref{qq-18}), (\ref{qq-28}), (\ref{qq-32}), we obtain
\begin{eqnarray}\label{qq-33}
\sup\limits_{t\in[t_0, \tau*)}|\nabla u(t)|^{2}+\int_{t_0}^{\tau*}(|\Delta u(t)|^{2}+|\nabla \partial_z u|^2)dt\leq C.
\end{eqnarray}
Taking the inner product of the equation (\ref{qq-2}) with $-\Delta \theta-\theta_{zz} $ in $H_2$,
similar to the above, we get
\begin{eqnarray*}
&&\frac{1}{2}\partial_{t}(|\nabla \theta|^{2}+|\theta_{z}|^{2}+ \gamma |\nabla \theta(z=0) |^{2}  )
+|\Delta \theta|^{2}+2( |\nabla \theta_{z} |^{2}+ \gamma|\nabla \theta(z=0) |^{2} )+|\theta_{zz}|^{2}\nonumber\\
&=&\int_{\mathcal{O}}(\alpha^{-1} u\cdot \nabla \theta+ \alpha^{-1} \Phi(u ) \theta_{z}-\beta Q)       (\Delta \theta+\theta_{zz} )\nonumber\\
&\leq& \varepsilon (| \Delta \theta|^{2}+|\theta_{zz}|^{2}+|\nabla \theta_{z}|^{2} )+C\beta^{2}|Q|^{2}+C \alpha^{-4}|\nabla u|^{2} |\Delta u|^{2} |\theta_{z}|^{2}\nonumber\\
&&+C (\alpha^{-2}|\tilde{u}|_{4}^{2}+\alpha^{-2}|\nabla \bar{u}|^{2}+\alpha^{-8}|\tilde{u}|_{4}^{8}+\alpha^{-8}|\nabla \bar{u}|^{8})|\nabla \theta|^{2}.
\end{eqnarray*}
Utilizing the Gronwall inequality, we get
\begin{eqnarray}\label{qq-40}
|\nabla \theta|^{2}+|\theta_{z}|^{2}+ \gamma |\nabla \theta(z=0) |^{2}
+\int^t_{t_0}[|\Delta \theta|^{2}+2( |\nabla \theta_{z} |^{2}+ \gamma|\nabla \theta(z=0) |^{2} )+|\theta_{zz}|^{2}]ds\leq  C.
\end{eqnarray}
In the following, we will complete our proof of the global well-posedness of stochastic PEs by three steps. Concretely,  we firstly prove the global existence of strong solution. Then, we show that the solution is continuous in the space $V$ with respect to $t$. At last, we prove the continuity in $V$ with respect to the initial data.  \\
\textbf{Step 1:}\quad the global existence of strong solution.\\
In the previous, we have obtained a priori estimates in $V.$ As we have indicated before that $[t_0, \tau_{*})$ is the maximal interval of existence of the solution of (\ref{qq-1})-(\ref{qq-6}),  we infer that $\tau_{*}=\infty,$ a.s..
Otherwise, if there exists $A\in \mathcal{F}$ such that $\mathbb{P}(A)>0$ and for fixed $\omega\in A, \tau_{*}(\omega)<\infty,$ it is clear that
\begin{eqnarray*}
\limsup\limits_{t\rightarrow \tau_{*}^{-}(\omega)}(\| u(t)\|_{1}+\|\theta(t)\|_{1} )=\infty,\ \ \mathrm{for}\ \mathrm{any}\ \omega\in A,
\end{eqnarray*}
which contradicts the priori estimates (\ref{qq-32}), (\ref{qq-33}) and (\ref{qq-40}).  Therefore $\tau_{*}=\infty,$ a.s.,  and the strong solution $(u, \theta)$ exists globally in time a.s..

\textbf{Step 2:} \quad  the continuity of strong solutions with respect to $t$.\\

Multiplying (\ref{qq-1}) by $\eta\in V_{1}$, integrating with respect to space variable, yields
\begin{eqnarray*}
\langle \partial_{t}A_{1}^{\frac{1}{2}} u, \eta \rangle= \langle \partial_{t}u, A_{1}^{\frac{1}{2}}  \eta \rangle&=&- \langle A_{1} u, A_{1}^{\frac{1}{2}}\eta\rangle - \alpha^{-1}\langle(u\cdot\nabla )u, A_{1}^{\frac{1}{2}}\eta\rangle\nonumber\\
&& -\alpha^{-1}\langle \Phi(u) \partial_{z}u, A_{1}^{\frac{1}{2}}  \eta \rangle-\langle  fu^{\bot}, A_{1}^{\frac{1}{2}}  \eta \rangle  \nonumber\\
&& + \alpha\beta^{-1}\langle\int_{-1}^{z}\nabla \theta dz',   A_{1}^{\frac{1}{2}}  \eta \rangle,
\end{eqnarray*}
where $ \langle \nabla p_{s},  A_{1}^{\frac{1}{2}}  \eta \rangle=0$ is used.
Taking a similar argument in (\ref{r-40}), we get
\begin{eqnarray*}
\langle \Phi(u) \partial_{z}u, A_{1}^{\frac{1}{2}}  \eta \rangle\leq \|u\|\|u\|_{2}|A_{1}^{\frac{1}{2}} \eta|.
\end{eqnarray*}
By the H\"{o}lder inequality and the Sobolev imbedding theorem, we have
\begin{eqnarray*}
\|\partial_{t}(A_{1}^{\frac{1}{2}}u)\|_{V_{1}'}\leq C(\|u\|_{2}+\|u\|\|u\|_{2}+|u|+|\nabla \theta| ).
\end{eqnarray*}
Since
\begin{eqnarray*}
u\in L^{\infty}([t_{0}, \mathcal{T}]; V_{1})\cap L^{2}([t_0, \mathcal{T}]; (H^{2}(\mathcal{O}))^{2}), \ \ \forall \mathcal{T}>t_{0},
\end{eqnarray*}
we obtain
\begin{eqnarray*}
A_{1}^{\frac{1}{2}}u\in L^{2}([t_0,\mathcal{T}]; V_{1} ),\quad  \partial_{t}(A_{1}^{\frac{1}{2}}u )\in L^{2}([t_0,\mathcal{T}]; V_{1}' ).
\end{eqnarray*}
Referring to Lemma \ref{l-1}, we deduce that
\[
A_{1}^{\frac{1}{2}}u\in C([t_0,\mathcal{T}]; H_{1})\ {\rm{or}}\  u\in C([t_0,\mathcal{T}]; V_{1})  \mathbb{P}-a.s..
\]
To study the regularity of $\theta,$ we choose $\xi\in V_{2}$. By (\ref{qq-2}) we have
\begin{eqnarray*}
\langle \partial_{t}A_{2}^{\frac{1}{2}} \theta ,      \xi \rangle& =&  \langle \partial_{t} \theta ,    A_{2}^{\frac{1}{2}}  \xi \rangle
=\langle  A_{2}\theta,   A_{2}^{\frac{1}{2}} \xi  \rangle- \alpha^{-1}\langle u\cdot \nabla \theta, A_{2}^{\frac{1}{2}}\xi     \rangle  \nonumber\\
&+&\alpha^{-1}\langle \Phi(u)\partial_{z}\theta,      A_{2}^{\frac{1}{2}}\xi   \rangle +\beta\langle Q,    A_{2}^{\frac{1}{2}}\xi \rangle.
\end{eqnarray*}
Taking a similar argument as above, we get
\begin{eqnarray*}
\langle \Phi(u)\partial_{z} \theta,      A_{2}^{\frac{1}{2}}\xi   \rangle \leq C\|u\|^{\frac{1}{2}}\|u\|_{2}^{\frac{1}{2}}
\|\theta\|^{\frac{1}{2}}\|\theta\|_{2}^{\frac{1}{2}}| A_{2}^{\frac{1}{2}}\xi |.
\end{eqnarray*}
Then by the H\"{o}lder inequality and the Sobolev imbedding theorem, we have
\begin{eqnarray*}
|\partial_{t} A_{2}^{\frac{1}{2}} \theta|_{V_{2}'}&\leq& C(|A_{2}\theta|+\alpha^{-1}\|u\|\|\theta\|_{2}\nonumber\\
&&+\alpha^{-1} \|u\|^{\frac{1}{2}}\|u\|_{2}^{\frac{1}{2}}
\|\theta\|^{\frac{1}{2}}\|\theta\|_{2}^{\frac{1}{2}} +\beta|Q|).
\end{eqnarray*}
In view of step one, we have
\begin{eqnarray*}
\theta\in L^{\infty}([t_0, \mathcal{T}]; V_{2})\cap L^{2}([t_0, \mathcal{T}]; H^{2}(\mathcal{O})), \ \ \forall \mathcal{T}>t_0.
\end{eqnarray*}
Therefore, by the same argument as above, we get
\begin{eqnarray*}
A_{2}^{\frac{1}{2}}\theta\in L^{2}([t_0,\mathcal{T}]; V_{2} ),\ \ \ \partial_{t}(A_{2}^{\frac{1}{2}}\theta )\in L^{2}([t_0,\mathcal{T}]; V_{2}' ).
\end{eqnarray*}
We deduce from Lemma \ref{l-1} that
\begin{eqnarray*}
A_{2}^{\frac{1}{2}}\theta\in C([t_0,\mathcal{T}]; H_{2})\  \mathrm{or}\  \theta\in C([t_0,\mathcal{T}]; V_{2}),\quad   \mathbb{P}-a.s..
\end{eqnarray*}
\textbf{Step 3:}\quad the continuity in $V$ with respect to the initial data.\\
Let $(\upsilon_{1}, T_{1} )$ and $(\upsilon_{2}, T_{2}) $ be two solutions of the system
(\ref{qq-1})-(\ref{qq-6}) with corresponding pressure $p_{b}{'}$ and $p_{b}{''},$ and initial data $((\upsilon_{0})_{1}, (T_{0})_{1})$ and
$((\upsilon_{0})_{2}, (T_{0})_{2}),$ respectively. Denote by $v=\upsilon_{1}-\upsilon_{2}, p_{b}=p_{b}{'}-p_{b}{''}$ and $T=T_{1}-T_{2}.$ Then
we have
\begin{eqnarray}\notag
&\partial_{t}v-\Delta v-\partial_{zz} v+\alpha^{-1} \upsilon_{1}\cdot \nabla v
+\alpha^{-1}(v\cdot \nabla) \upsilon_{2}+\alpha^{-1}\Phi( \upsilon_{1})v_{z}+\alpha^{-1}\Phi(v)\partial_{ z}\upsilon_{2}&\\
\label{r-1}
& +fk\times v+\alpha\nabla p_{b}
-\alpha\beta^{-1}\int_{-1}^{z}\nabla Tdz'=0,&\\
\label{r-2}
&\partial_{t}T-\Delta T-\partial_{zz}T+\alpha^{-1}\upsilon_{1}\cdot \nabla  T
+\alpha^{-1}(v\cdot \nabla) T_{2}
+\alpha^{-1}\Phi (\upsilon_{1} )T_{z}+\alpha^{-1}\Phi (v) \partial_{z}T_{2} =0,&\\
&\int_{-1}^{0}\nabla\cdot vdz=0,&\\
\label{r-3}
&v(x,y,z,t_0)=(\upsilon_{0})_{1}-(\upsilon_{0})_{2},\ \ \ \ T(x,y,z,t_0)=(T_{0})_{1}-(T_{0})_{2},&\\
\label{r-4}
&(v,T)\ {\rm{satisfies\ the \ boundary\ value\ conditions }}\ (\ref{qq-4})-(\ref{qq-6})&.
\end{eqnarray}
Multiplying $L_{1}v$ in equation (\ref{r-1}) and integrating with respect to spatial variable, we have
\begin{eqnarray}
&&\frac{1}{2}\partial_{t} (|\nabla v|^{2}+|\partial_{z} v|^{2})+| \Delta v|^{2}+|\partial_{z} v|^{2}+|\nabla v_{z}|^{2}\nonumber\\
&=&-\alpha^{-1}\int_{\mathcal{O}}(\upsilon_{1}\cdot \nabla v  )\cdot L_{1}v -\alpha^{-1}\int_{\mathcal{O}}\Phi(\upsilon_{1}  ) v_{z}\cdot L_{1}v\nonumber\\
&&\ -\alpha^{-1}\int_{\mathcal{O}}(\Phi(v)\partial_{z} \upsilon_{2})\cdot L_{1}v -\alpha^{-1}\int_{\mathcal{O}}[(v\cdot\nabla )\upsilon_{2}]\cdot L_{1}v  \nonumber\\
&& \ -\int_{\mathcal{O}}(fk\times v)\cdot L_{1}v
+\alpha \beta^{-1}\int_{\mathcal{O}}(\int_{-1}^{z}\nabla T dz' )\cdot L_{1}v\nonumber\\
&:=& K_1(t)+K_2(t)+K_3(t)+K_4(t)+K_5(t)+K_6(t).
\end{eqnarray}
Using the Agmon inequality and the H\"{o}lder inequality, we obtain
\begin{eqnarray*}
K_{1}(t)&\leq&\alpha^{-1} |v_{1}|_{\infty}|\nabla v||L_{1}v|\nonumber\\
&\leq& C\alpha^{-2}\|v_{1}\|^{\frac{1}{2}}\|v_{1}\|_{2}^{\frac{1}{2}}\|v\| \|v\|_{2}\nonumber\\
&\leq& C\alpha^{-2}\|v_{1}\|\|v_{1}\|_{2}\|v\|^{2}+\varepsilon \|v\|_{2}^{2}.
\end{eqnarray*}
Applying the H\"{o}lder inequality, the interpolation inequality and the Sobolev embedding theorem,  we get
\begin{eqnarray*}
K_{2}(t)&\leq& \alpha^{-1}\int_{\mathcal{O}}|\int_{-1}^{0} \nabla \cdot v_{1}dz |\cdot |v_{z}|\cdot|L_{1}v |\nonumber\\
&\leq & \alpha^{-1}|\nabla \cdot \bar{v}_{1} |_{L^{4}(M)}\int_{-1}^{0}|v_{z}|_{L^{4}(M)}|L_{1}v|_{L^{2}(M)}dz\nonumber\\
&\leq & C\alpha^{-1}\|v_{1}\|^{\frac{1}{2}}\|v_{1}\|_{2}^{\frac{1}{2}}\|v\|^{\frac{1}{2}}
\|v\|_{2}^{\frac{1}{2}}\|v\|_{2}\nonumber\\
&\leq &\varepsilon \|v\|_{2}^{2}+C\alpha^{-2}\|v_{1}\|^{2}\|v_{1}\|_{2}^{2}\|v\|^{2}.
\end{eqnarray*}
Similar to the above, we obtain
\begin{eqnarray*}
K_{3}(t)&\leq&\alpha^{-1} \int_{\mathcal{O}} |\int_{-1}^{0}\nabla\cdot v dz|\cdot|\partial_{z} v_{2} |\cdot |L_{1}v|\nonumber\\
&\leq&\alpha^{-1}|\nabla \cdot \bar{v} |_{L^{4}(M)}\int_{-1}^{0}|\partial_{z}v_{2}|_{L^{4}(M)}|L_{1}v|_{L^{2}(M)}dz\nonumber\\
&\leq&C\alpha^{-1}\|v\|^{\frac{1}{2}}\|v\|_{2}^{\frac{1}{2}}\int_{-1}^{0}\|v_{2}\|_{H^{1}(M)}^{\frac{1}{2}}\|v_{2}\|_{H^{2}(M)}^{\frac{1}{2}}\|v\|_{2}dz\nonumber\\
&\leq&\varepsilon \|v\|_{2}^{2}+C\alpha^{-4}\|v\|^{2}\|v_{2}\|^{2}\|v_{2}\|_{2}^{2}.
\end{eqnarray*}
With the help of the H\"{o}lder inequality, we deduce that
\begin{eqnarray*}
K_{4}(t)&\leq& \alpha^{-1}|v|_{4}|\nabla v_{2} |_{4}|L_{1}v|_{2}\leq \varepsilon \|v\|_{2}^{2}+C\alpha^{-2}\|v\|^{2}\|v_{2}\|_{2}^{2},
\end{eqnarray*}
and
\begin{eqnarray*}
K_{5}(t)+K_{6}(t)\leq \varepsilon  \|v\|_{2}^{2}+ C|v|_{2}^{2}+C\alpha^{2}\beta^{-2}\|T\|^{2}.
\end{eqnarray*}
By the boundary conditions, we know that
 $|T|^{2}$ is small than $|T_{z}|^{2}+|T(z=0)|_{L^2(M)}^{2}$, then
$\|T\|^{2}$ is equivalent to $|\nabla T|^{2}+|T_{z}|^{2}+|T(z=0)|_{L^2(M)}^{2}.$ Keeping this in mind and taking an inner product of the equation (\ref{r-2}) with $L_{2}T$, we have
\begin{eqnarray*}
&&\frac{1}{2}\partial_{t}( |\nabla T|^{2}+|T_{z}|^{2}+\gamma|T(z=0)|^{2} )
+ |\Delta T |^{2}+|T_{zz}|^{2}+|\nabla T_{z}|^{2}+\gamma|\nabla T(z=0) |^{2}\nonumber\\
&=&-\alpha^{-1}\int_{\mathcal{O}}(\upsilon_{1}\cdot \nabla  T)L_{2}T
-\alpha^{-1}\int_{\mathcal{O}}(v\cdot \nabla T_{2})L_{2}T\nonumber\\
 &&-\alpha^{-1}\int_{\mathcal{O}}\Phi (\upsilon_{1} )T_{z}L_{2}T+\alpha^{-1}\Phi (v)\partial_{z}T_{2}L_{2}T\\
 &:=&J_{1}(t)+J_2(t)+J_3(t)+J_4(t).
\end{eqnarray*}
By the Agmon inequality, we get that
\begin{eqnarray*}
J_{1}(t)+J_{2}(t)&\leq& \alpha^{-1}|v_{1}|_{\infty}|\nabla T||L_{2}T|+\alpha^{-1}|v|_{4}|\nabla T_{2}|_{4}|L_{2}T|\nonumber\\
&\leq&C\alpha^{-1}\|v_{1}\|^{\frac{1}{2}}\|v_{1}\|^{\frac{1}{2}}_2\|T\|\|T\|_{2}+C \alpha^{-1} \|v\|\|T_{2}\|_{2}\|T\|_{2}\nonumber\\
&\leq& \varepsilon \|T\|_{2}^{2}+C\alpha^{-2}\|v\|^{2}\|T_{2}\|_{2}^{2}+C\alpha^{-2}\|v_{1}\|\|v_{1}\|_{2}\|T\|^{2}.
\end{eqnarray*}
Taking an similar argument as $K_{2}(t)$, we obtain
\begin{eqnarray*}
J_{3}(t)&\leq&\alpha^{-1}\int_{\mathcal{O}} \Big(\int_{-1}^{0}|\nabla \cdot v_{1}|dz\Big)\cdot |T_{z}|\cdot |L_{2}T| \nonumber\\
&\leq& \alpha^{-1}\Big{(}\int_{-1}^{0}|L_{2}T|_{L^{2}(M)}|T_{z}|_{L^{4}(M)}dz \Big{)}|\int_{-1}^{0}|\nabla \cdot v_{1} |dz|_{L^{4}(M)}\nonumber\\
&\leq&C\alpha^{-1}\int_{-1}^{0}\|T\|_{H^{2}(M)}^{\frac{3}{2}}\|T\|_{H^{1}(M)}^{\frac{1}{2}}dz\int_{-1}^{0}|\nabla \cdot v_{1}|_{L^{4}(M)}dz\nonumber\\
&\leq&C\alpha^{-1}\|T\|_{2}^{\frac{3}{2}}\|T\|^{\frac{1}{2}}\int_{-1}^{0}\|v_{1}\|_{H^{1}(M)}^{\frac{1}{2}}\|v_{1}\|_{H^{2}(M)}^{\frac{1}{2}}dz\nonumber\\
&\leq&C\alpha^{-1}\|T\|_{2}^{\frac{3}{2}}\|T\|^{\frac{1}{2}}\|v_{1}\|^{\frac{1}{2}}\|v_{1}\|_{2}^{\frac{1}{2}}\nonumber\\
&\leq&\varepsilon \|T\|_{2}^{2}+C\alpha^{-4}\|T\|^{2}\|v_{1}\|^{2}\|v_{1}\|_{2}^{2}.
\end{eqnarray*}
Similarly, we can prove that
\begin{eqnarray*}
J_{4}(t)&\leq&\alpha^{-1} \int_{\mathcal{O}}(\int_{-1}^{0}|\nabla \cdot v|)|\partial_{z}T_{2}  |\cdot|L_{2}T| \\
&\leq&\alpha^{-1}\Big{(} \int_{-1}^{0} |L_{2}T|_{L^{2}(M)}|\partial_{z}T_{2} |_{L^{4}(M)}dz\Big{)}|\int_{-1}^{0}|\nabla \cdot v|dz |_{L^{4}(M)}\nonumber\\
&\leq&\alpha^{-1}\Big{(} \int_{-1}^{0}\|T\|_{H^{2}(M)}\|T_{2}\|_{H^{1}(M)}^{\frac{1}{2}}\|T_{2}\|_{H^{2}(M)}^{\frac{1}{2}}dz\Big{)}\int_{-1}^{0}|\nabla \cdot v|_{L^{4}(M)}dz\nonumber\\
&\leq& C\alpha^{-1}\|T\|_{2}\|T_{2}\|^{\frac{1}{2}}\|T_{2}\|_{2}^{\frac{1}{2}}\|v\|^{\frac{1}{2}}\|v\|_{2}^{\frac{1}{2}}\nonumber\\
&\leq& \varepsilon \|T\|_{2}^{2}+\varepsilon \|v\|_{2}^{2}+C\alpha^{-4}\|T_{2}\|^{2}\|T_{2}\|_{2}^{2}\|v\|^{2}.
\end{eqnarray*}
Denote
\begin{eqnarray*}
\eta(t)&:=& \|v(t)\|^{2}+\|T(t)\|^{2},\\
\xi(t)&:=&\alpha^{-2}\|v_{1}\|\|v_{1}\|_{2}+\alpha^{-4}\|v_{1}\|^{2}\|v_{1}\|_{2}^{2}\\
&&+\alpha^{-4} \|v_{2}\|^{2}\|v_{2}\|_{2}^{2}+\alpha^{2} \|v_{2}\|_{2}^{2}+\alpha^{2}\|T_{2}\|_{2}^{2}+\alpha^{-4}\|T_{2}\|^{2} \|T_{2}\|_{2}^{2}  +1.
\end{eqnarray*}
Notice that $|\nabla v|^{2}+|\partial_{z} v|^{2}$ is equivalent to $\|v\|^{2}$ and $|\nabla T|^{2}+|T_{z}|^{2}+|T(z=0)|_{L^2(M)}^{2}$ is equivalent to  $\|T\|^{2}$, letting $\varepsilon$ be small enough, we deduce from the above estimates of $K_{1}-K_{6}$ and $J_{1}-J_{4}$ that
\begin{eqnarray}
\frac{d \eta(t)}{dt}+\|v\|_{2}^{2}+\|T\|_{2}^{2}\leq \eta(t)\xi(t).
\end{eqnarray}
Since $(v_{i}(t),T_{i}(t)),i=1,2,$ is the solution of stochastic PEs in sense of Definition \ref{d-1}, we have
\begin{eqnarray*}
\int_{t_0}^{t}\xi(s)ds< \infty,\ \ a.s.,\ \mathrm{for}\ \mathrm{all}\ t\in(t_0, \infty),
\end{eqnarray*}
which implies that
\begin{eqnarray*}
\eta(t)\leq \eta(t_0)e^{C\int_{t_0}^{t}\xi(s)ds}.
\end{eqnarray*}
Therefore, we proved that for any $t\in(t_0, \infty) , (u(t), \theta(t))$ is Lipschitz continuous in $V_{1}\times V_{2}$ with respect to
the initial data $(u(t_0), \theta(t_0)),$ which  is equivalent to that strong solution $ (\upsilon(t), T(t))$ of (\ref{e-1})-(\ref{e-5}) is Lipschitz continuous in $V_{1}\times V_{2}$ with respect to the initial data $(v_{t_0}, T_{t_0}),$ for any $t\in (t_0, \infty).$
\hspace{\fill}$\square$

\end{flushleft}
\begin{Rem}\ \\
 \textbf{(1)}\quad
 In the above theorem, we have obtained the continuity of the strong solution with respect to initial data in $(H^{1}(\mathcal{O}))^{3}$. This is a key to prove the compact property of the solution operator in $V$. Notice that the authors only proved the strong solution is Lipschitz continuous in the space $( L^{2}(\mathcal{O}))^{3}$ with respect to the initial data in \cite{GH}, which is not enough to obtain the asymptotical behavior in $(H^{1}(\mathcal{O}))^{3}$.\\
\textbf{(2)}\quad With the help of Lemma \ref{l-1}, we have established the the continuity of strong solution with respect to time in $V$ and a priori estimates to prove the compact property of the solution operator in $V$.\\
\textbf{(3)}\quad We release the regularity of $Q$ from $H^{1}(\mathcal{O})$ to $L^{2}(\mathcal{O})$, which is more natural.
\end{Rem}

\section{Existence of random attractor}
In this section, we establish the existence of random attractor. Firstly, we recall some preliminaries from \cite{CDF}.

Denote by $C_{0}(\mathbb{R}; X )$ the space of continuous functions with values in $X$ and equal to $0$ at $t=0$.
Let $(X, d)$ be a polish space and $(\tilde{\Omega}, \tilde{\mathcal{F}}, \tilde{\mathbb{P}}  )$ be a probability space, where $ \tilde{\Omega}$ is the two-sided Wiener space $C_{0}(\mathbb{R}; X )$.
\begin{definition}\label{d-3}
A family of maps
$S(t,s;\omega):X\rightarrow X,\ \ -\infty<s \leq t< \infty$, parametrized by $\omega\in \tilde{\Omega},$ is said to be a stochastic flow, if $ \tilde{\mathbb{P}}$-a.s.,
\begin{description}
  \item[(i)] $S(t,r;\omega)S(r,s;\omega)x= S(t,s;\omega)x$\ for\ all\ $ s\leq r\leq t,\  x\in X,$
  \item[(ii)]  $\ \ S(t,s;\omega)$\ is\ continuous\ in\ $X,$ for\ all\ $s\leq t,$
  \item[(iii)]for all $s<t$ and $x\in X$, the mapping
\[
\omega\mapsto S(t,s;\omega)x
\]
is measurable from $(\tilde{\Omega},\tilde{\mathcal{F}})$ to $(X, \mathcal{B}(X )  )$ where $\mathcal{B}(X ) $ is the Borel $\sigma$-algebra of $ X$,
  \item[(iv)] for all $t, x\in X,$ the mapping $s\mapsto S(t,s;\omega)$ is right continuous at any point.
\end{description}
\end{definition}
\begin{definition}
A set-valued map $K: \tilde{\Omega}\rightarrow 2^{X}$ taking values in the closed subsets of $X$ is said to be measurable if for each $x\in X$ the map
$\omega\mapsto d(x, K(\omega))$ is measurable, where
\[
d(A,B)=\sup\{\inf\{d(x,y):y\in B \}:x\in A \} \ {\rm{for}} A, B \in 2^{X}, A,B\neq \emptyset,
\]
and $d(x,B)=d(\{x\},B).$ Since $ d(A,B)=0$ if and only if $A\subset B$, $d$ is not a metric.
\end{definition}
\begin{definition}
 A closed set valued measurable map $K:\tilde{\Omega}\rightarrow 2^{X}$ is called a random closed set.
\end{definition}
\begin{definition}
Given $t\in \mathbb{R}$ and $\omega\in \tilde{\Omega}, K(t,\omega)\subset X$ is called an attracting set at time $t$ if for all bounded sets $B\subset X,$
\[
d(S(t,s;\omega)B, K(t,\omega) )\rightarrow 0,\ \ provided\ s\rightarrow -\infty.
\]
Moreover, if for all bounded sets $B\subset X,$ there exists $t_{B}(\omega)$ such that for all $s\leq t_{B}(\omega)$
\[
S(t,s;\omega)B\subset K(t,\omega),
\]
we say $K(t,\omega) $ is an absorbing set at time $t.$
\end{definition}
Let $\{\vartheta_{t}:\tilde{\Omega}\rightarrow \tilde{\Omega}   \}, t\in T=\mathbb{R},$ be a family of measure preserving transformations of the probability space $(\tilde{\Omega}, \tilde{\mathcal{F}},\tilde{ \mathbb{P}} )$ such that for all $s< t$ and $\omega\in \tilde{\Omega}$, the following
\begin{description}
  \item[(a)] $(t,\omega)\rightarrow \vartheta_{t}\omega$ is measurable,
  \item[(b)] $\vartheta_{t}(\omega)(s)=\omega(t+s)-\omega(t)$,
  \item[(c)] $S(t,s;\omega)x=S(t-s,0;\vartheta_{s}\omega)x$,
\end{description}
hold. Then, $(\vartheta_{t} )_{t\in T}$ is a flow and
$((\tilde{\Omega}, \tilde{\mathcal{F}},\tilde{ \mathbb{P}} ), (\vartheta_{t} )_{t\in T} )$ is a measurable dynamical system.

\begin{definition}
Given a bounded set $B\subset X$, the set
\begin{eqnarray*}
\Omega(B,t,\omega)=\bigcap\limits_{T\leq t}\overline{\bigcup\limits_{s\leq T}S(t, s,\omega)B}
\end{eqnarray*}
is said to be the $\Omega$-limit set of $B$ at time $t$. Obviously, if denote $\Omega(B,0,\omega)=\Omega(B,\omega),$ we have
$\Omega(B,t,\omega)=\Omega(B,\vartheta_{t}\omega).$
\end{definition}
It's easy to identify
\begin{eqnarray*}
\Omega(B,t,\omega)=\Big\{x\in X: {\rm{there\ exists}}\  s_{n}\rightarrow -\infty\ {\rm{and}}\ x_{n}\in B\ {\rm{such\ that}}\ \lim\limits_{n\rightarrow\infty}S(t,s_{n},\omega)x_{n}=x\Big\}.
\end{eqnarray*}
Furthermore, if there exists a compact attracting set $K(t,\omega)$ at time $t,$ it is not difficult to check that $\Omega(B,t,\omega)$ is a nonempty compact subset of $X$ and $\Omega(B,t,\omega)\subset K(t,\omega). $
\begin{definition}\label{d-2}
For all $t\in \mathbb{R}$ and $\omega\in \tilde{\Omega},$  a random closed set $\omega\rightarrow \mathcal{A}(t,\omega)$ is called the random attractor, if $\tilde{\mathbb{P}}-$a.s.,
\begin{description}
  \item[(1)] $\mathcal{A}(t,\omega)$ is a nonempty compact subset of $X,$
  \item[(2)] $\mathcal{A}(t,\omega)$ is the minimal closed attracting set,
i.e., if $\tilde{\mathcal{A}}(t,\omega)$ is another closed attracting set, then $\mathcal{A}(t,\omega)\subset \tilde{\mathcal{A}}(t,\omega),$
  \item[(3)] it is invariant,  in the sense that, for all $s\leq t,$ $
S (t,s;\omega)\mathcal{A}(s,\omega)=\mathcal{A}(t,\omega).
 $
\end{description}
\end{definition}
Let
$\mathcal{A}(\omega)=\mathcal{A}(0,\omega)$,
then the invariance property can be written as
$$S(t,s;\omega)\mathcal{A}(\vartheta_{s}\omega)=\mathcal{A }(\vartheta_{t}\omega).$$
We will prove the existence of the random attractor using Theorem 2.2 in $\cite{CDF} $.
For the convenience of reference, we cite it here.
\begin{theorem}\label{thm-1}
Let $(S(t, s; \omega))_{t\geq s, \omega\in \tilde{\Omega}}$ be a stochastic dynamical system
satisfying \textbf{(i)-(iv)}. Assume that there exists a family of measure preserving mappings $ \vartheta_{t}, t\in \mathbb{R}$ such that\textbf{(a)}-  \textbf{(c)} hold and there exists a compact attracting set $K(\omega)$ at time $0$, for $\tilde{\mathbb{P}}$-a.s..  Set
$$\mathcal{A}(\omega)=\overline{\bigcup_{B\subset X,\\ B \ bounded} \Omega(B,\omega) } $$
where the union is taken over all the bounded subsets of $X$. Then we have  $\tilde{\mathbb{P}}  $-a.s.,
\begin{description}
  \item[(1)] $\mathcal{A}(\omega)$ is a nonempty compact subset of $X$. If $X$ is connected,
it is a connected subset of $K(\omega)$.
  \item[(2)] The family $\mathcal{A}(\omega),\ \omega\in \Omega$, is measurable.
  \item[(3)] $\mathcal{A}(\omega)$ is invariant in the sense that
$S(t,s;\omega)\mathcal{A}(\vartheta_{s}\omega)= \mathcal{A}(\vartheta_{t}\omega),\ \ s\leq t.$
  \item[(4)] It attracts all bounded sets from $-\infty$: for bounded $B\subset X$ and $\omega\in \tilde{\Omega}$
\begin{eqnarray*}
d(S(t,s;\omega)B, \mathcal{A}(\vartheta_{t}\omega))\rightarrow 0,\ \ when\ s\rightarrow -\infty.
\end{eqnarray*}
Moreover, it is the minimal closed set with this property: if $\tilde{\mathcal{A}}(\vartheta_{t}\omega)$ is a closed attracting set, then $\mathcal{A}(\vartheta_{t}\omega)\subset \tilde{\mathcal{A}}(\vartheta_{t}\omega).$
  \item[(5)]  For any bounded set $ B\subset X,\ d(S(t,s;\omega)B, \mathcal{A}(\vartheta_{t}\omega))\rightarrow 0$ in probability when $t\rightarrow \infty.$
\end{description}
And if the time shift $\vartheta_{t},t\in \mathbb{R}$ is ergodic,
\begin{description}

  \item[(6)] there exists a bounded set $B\subset X$ such that
\begin{eqnarray*}
\mathcal{A}(\omega)= \mathcal{A}(B, \omega),
\end{eqnarray*}
  \item[(7)] $ \mathcal{A}(\omega)$ is the largest compact measurable set which is invariant in sense of Definition \ref{d-2}.
\end{description}
\end{theorem}

Before showing the existence of random attractor, we recall the Aubin-Lions Lemma, which is vital to the proof of Theorem \ref{thm-3}.
\begin{lemma}\label{lem-1}
Let $B_{0}, B, B_{1}$ be Banach spaces such that $B_{0}, B_{1}$ are reflexive and $B_{0}\overset{c}{\subset}B\subset B_{1}.$
For $0<T<\infty,$ set
\begin{eqnarray*}
X:=\Big{\{}h \Big{|}h\in L^{2}([0,T]; B_{0}), \frac{dh}{dt}\in L^{2}([0,T]; B_{1})\Big{\}}.
\end{eqnarray*}
Then $X$ is a Banach space equipped with the norm $|h|_{L^{2}([0,T]; B_{0})}+|h'|_{L^{2}([0,T]; B_{1})}.$ Moreover, $$X \overset{c}\subset {L^{2}([0,T]; B)}.$$
\end{lemma}
The main result is:
\begin{theorem}\label{thm-3}
Let $Q\in L^{2}(\mathcal{O}), \upsilon_{0}\in V_{1}, T_{0}\in V_{2} $. Then the solution operator $(S(t,s;\omega))_{t\geq s,\omega\in \tilde{\Omega}} $ of 3D stochastic PEs (\ref{e-1})-(\ref{e-5}): $S(t,s;\omega)(\upsilon_{s}, T_{s})=(\upsilon(t), T(t) ) $ satisfies \textbf{(i)-(iv)} in Definition \ref{d-3} and possesses a compact absorbing ball $\mathcal{B}(0,\omega)$ in $V$ at time $0$.  Furthermore, for $\tilde{\mathbb{P}}$-a.s. $\omega\in \Omega$, set
$$\mathcal{A}(\omega)=\overline{\bigcup_{B\subset V}\Omega(B,\omega) } $$
where the union is taken over all the bounded subsets of $V$. Then $\mathcal{A}(\omega)$ is the random attractor of stochastic PEs (\ref{e-1})-(\ref{e-5}) and possesses the properties \textbf{(1)-(7)} of Theorem \ref{thm-1} with space $X$ replaced by space $V.$
\end{theorem}

\begin{flushleft}
\textbf{Proof.}\quad
Denote by $w=(w_{1}:=\sum_{k=1}^{n}\alpha_{k}w_{k}^{1}, w_{2}:=\sum_{k=1}^{n}\beta_{k}w_{k}^{2})$ the $\mathbb{R}^{2}-$valued Brownian motion, which has a version $\omega$ in $ C_{0}(\mathbb{R}, \mathbb{R}^{2}):=\tilde{\Omega}$, the space of continuous functions which are $\mathrm{zero}$ at $\mathrm{zero}.$ In the following, we consider a canonical version of $w$ given by the probability space $(C_{0}(\mathbb{R},\mathbb{R}^{2}), \mathcal{B}(C_{0}(\mathbb{R},  \mathbb{R}^{2})), \tilde{\mathbb{P}} )$, where $\tilde{\mathbb{P}}$ is the Wiener-measure generated by $w$ and $\mathcal{B}(C_{0}(\mathbb{R},  \mathbb{R}^{2})) $ is the family of Borel subsets of $ C_{0}(\mathbb{R},  \mathbb{R}^{2})$.
Now, define the stochastic flow $(S(t,s;\omega))_{t\geq s,\ \omega\in \tilde{\Omega}}$ by
\begin{eqnarray}
S(t,s;\omega)(\upsilon_{s},T_{s})=(\alpha^{-1}(t)u(t,\omega_{1}), \beta^{-1}(t)\theta(t, \omega_{2})),
\end{eqnarray}
where $(\upsilon, T)$ is the strong solution to  (\ref{e-1})-(\ref{e-5}) with $(\upsilon_{s}, T_{s})=(\alpha^{-1}(s)u_{s}(s,\omega_{1}),\beta^{-1}(s) \theta_{s}(s,\omega_{2}  ))$ and $(u, \theta)$ is the strong solution to (\ref{qq-1})-(\ref{qq-6}).
It can be checked that assumptions \textbf{(i)-(iv)} and \textbf{(a)-(c)} of stochastic dynamics are satisfied with $X=V$. Indeed, properties \textbf{(i),(ii),(iv)} of the solution operator $(S(t,s;\omega))_{t\geq s,\omega\in \tilde{\Omega}} $ follows by Theorem \ref{thm-2} and property \textbf{(iii)} of the solution operator also holds with the help of the global existence of strong solution to  (\ref{e-1})-(\ref{e-5}) rest upon Faedo-Galerkin method.
Furthermore,
$(\tilde{\Omega}, \mathcal{B}(C_{0}(\mathbb{R}, \mathbb{R}^{2} )), \tilde{\mathbb{P}}, \vartheta )$ is an ergodic metric dynamical system.
\par
In the following, we will prove the existence of the random attractor.
Let $(u(t,\omega;t_{0},u_{0}), \theta(t,\omega;t_{0},\theta_{0})) $  be the solution to $(\ref{qq-1})-(\ref{qq-6})$ with  initial value $u(t_{0})=u_{0}$ and $\theta(t_{0})=\theta_{0}.$
By the law of the iterated logarithm, we have
\begin{eqnarray}
 \lim\limits_{t\rightarrow -\infty}\frac{\sum_{k=1}^{n}\alpha_{k}w_{k}^{1}}{t}=\lim\limits_{t\rightarrow -\infty}\frac{\sum_{k=1}^{n}\beta_{k}w_{k}^{2}}{t}=0.
\end{eqnarray}
Obviously, $t\rightarrow \beta^{2}(t)e^{\lambda t}$ is pathwise integrable over $(-\infty, 0]$, where $\lambda$ is positive. And we have
\begin{eqnarray}\label{r-8}
 \lim\limits_{t\rightarrow -\infty} \beta^{2}(t)e^{\lambda t}=0,\ \ \ \tilde{\mathbb{P}}-a.e..
\end{eqnarray}
In view of (\ref{qq-8}), (\ref{qq-9}) and (\ref{r-8}), for $\tilde{\mathbb{P}}-$a.e. $\omega\in \tilde{\Omega} $, there exists a random variable $r_{1}(\omega)$, depending only on
$\lambda$ such that for arbitrary $\rho >0$ there exists $t(\omega)\leq -4$ such that for all $t_{0}\leq t(\omega)$ and $(u_{0}, \theta_{0})\in V$ with $\|u_{0}\|+\|\theta_{0}\| \leq \rho,$
$ \theta(t,\omega;t_{0},\theta_{0}) $ satisfies
\begin{eqnarray}\label{r-9}
\sup\limits_{t\in [-4, 0]}|\theta(t,\omega;t_{0},\theta_{0})|^{2}+\int_{-4}^{0}\|\theta(s)\|^{2}ds \leq r_{1}(\omega).
\end{eqnarray}
In view of (\ref{r-10}) and (\ref{r-9}), taking a similar argument as (\ref{r-9}), for $\tilde{\mathbb{P}}-$a.e. $\omega\in \tilde{\Omega} $, we deduce that
 there exists random variable $r_{2}(\omega)$, depending only on
$\lambda$ such that for arbitrary $\rho >0$ there exists $t(\omega)\leq -4$ such that for all $t_{0}\leq t(\omega)$ and $(u_{0}, \theta_{0})\in V$ with $\|u_{0}\|+\|\theta_{0}\| \leq \rho,$
$ u(t,\omega;t_{0},u_{0}) $ satisfies
\begin{eqnarray}\label{r-16}
\sup\limits_{t\in [-4, 0]}|u(t,\omega;t_{0},u_{0})|^{2}+\int_{-4}^{0}\|u(s)\|^{2}ds \leq r_{2}(\omega).
\end{eqnarray}

By (\ref{e-12}), repeating the argument as in (\ref{r-9}), for $\tilde{\mathbb{P}}-$a.e. $\omega\in \tilde{\Omega} $, there exists random variable $r_{3}(\omega)$, depending only on
$\lambda$ such that for arbitrary $\rho >0$ there exists $t(\omega)\leq -4$ such that for all $t_{0}\leq t(\omega)$ and $(u_{0}, \theta_{0})\in V$ with $\|u_{0}\|+\|\theta_{0}\| \leq \rho,$
$ \theta(t,\omega;t_{0},\theta_{0}) $ satisfies
\begin{eqnarray}\label{r-17}
\sup\limits_{t\in [-4, 0]}|\theta(t,\omega;t_{0},\theta_{0})|_{4}^{2}\leq r_{3}(\omega).
\end{eqnarray}
Integrating  (\ref{r-12}) with respect to time over $[t,-3]$ yields,
\begin{eqnarray}\label{r-13}
|\tilde{u}(-3)|_{4}^{2}\leq \Big( |\tilde{u}(t)|_{4}^{2}+C\int_{t}^{-3}\alpha^{2}(s)\beta^{-2}(s)|\theta(s)|_{4}^{2}ds\Big)
{\rm{e}}^{C\int_{t}^{-3}\alpha^{-2}(s)\|u(s)\|^{2}(1+\alpha^{-2}(s)|u(s)|^{2})ds},
\end{eqnarray}
Integrating  (\ref{r-13}) with respect to $t$ over $[-4,-3]$, we obtain
\begin{eqnarray}
|\tilde{u}(-3)|_{4}^{2}&\leq& \Big( \int_{-4}^{-3}|\tilde{u}(t)|_{4}^{2}ds +C\int_{-4}^{-3}\alpha^{2}(s)\beta^{-2}(s)|\theta(s)|_{4}^{2}ds\Big)\nonumber\\
&& \times {\rm{e}}^{C\int_{-4}^{-3}\alpha^{-2}(s)\|u(s)\|^{2}(1+\alpha^{-2}(s)|u(s)|^{2})ds}.
\end{eqnarray}
Therefore, by virtue of (\ref{r-9})-(\ref{r-17}), we conclude that for $\tilde{\mathbb{P}}-$a.e. $\omega\in \tilde{\Omega} $, there exists random variable $C_{1}(\omega)$, depending only on
$\lambda$ such that for arbitrary $\rho >0$ there exists $t(\omega)\leq -4$ such that for all $t_{0}\leq t(\omega)$ and $(u_{0}, \theta_{0})\in V$ with $\|u_{0}\|+\|\theta_{0}\| \leq \rho,$
$\tilde {u}(-3,\omega;t_{0},u_{0}) $ satisfies
\begin{eqnarray}\label{r-18}
|\tilde{u}(-3)|_{4}^{2}\leq C_{1}(\omega).
\end{eqnarray}
Integrating (\ref{r-12}) with respect to time over $[-3,  t]$, we get
\begin{eqnarray}
|\tilde{u}(t)|_{4}^{2}\leq \Big( |\tilde{u}(-3)|_{4}^{2}+C\int_{-3}^{t}\alpha^{2}(s)\beta^{-2}(s)|\theta(s)|_{4}^{2}ds\Big)
{\rm{e}}^{C\int_{-3}^{t}\alpha^{-2}(s)\|u(s)\|^{2}(1+\alpha^{-2}(s)|u(s)|^{2})ds}.
\end{eqnarray}

In view of (\ref{r-9})-(\ref{r-17}) and (\ref{r-18}), for $\tilde{\mathbb{P}}-$a.e. $\omega\in \tilde{\Omega} $, we deduce that
 there exists random variable $r_{4}(\omega)$  such that for arbitrary $\rho >0$ there exists $t(\omega)\leq -3$ such that for all $t_{0}\leq t(\omega)$ and $(u_{0}, \theta_{0})\in V$ with $\|u_{0}\|+\|\theta_{0}\| \leq \rho,$
$ u(t,\omega;t_{0},u_{0}) $ satisfies
\begin{eqnarray}\label{r-24}
\sup\limits_{t\in [-3, 0]}|\tilde{u}(t,\omega;t_{0},\tilde{u}_{0})|_{4}^{2}\leq r_{4}(\omega).
\end{eqnarray}

Taking integration of (\ref{r-20}) with respect to time over $[-3, 0]$ yields,
\begin{eqnarray}\notag
&&\int_{-3}^{0}(|\nabla (|\tilde{u}|^{2})|^{2}+|\partial_{z}(|\tilde{u}|^{2})|^{2}+|\tilde{u}|\nabla \tilde {u}||^{2}
+|\tilde{u}|\partial_{z}\tilde{u}||^{2})dt\\
\label{r-22}
&\leq& |\tilde{u}(-3)|_{4}^{2}+C\int_{-3}^{0}\alpha^{2}\beta^{-2}|\theta|_{4}^{2}|\tilde u|_{4}^{2}+C\int_{-3}^{0}(\alpha^{-2}\|u\|^{2}+\alpha^{-4}|u|^{2}\|u\|^{2} )|\tilde{u}|_{4}^{2}.
\end{eqnarray}

By (\ref{r-16}), (\ref{r-17}) and (\ref{r-22}), for $\tilde{\mathbb{P}}-$a.e. $\omega\in \tilde{\Omega} $, we infer that
that there exists random variable $C_{2}(\omega)$
such that for $\rho >0$ there exists $t(\omega)\leq -3$ such that for all $t_{0}\leq t(\omega)$ and $(u_{0}, \theta_{0})\in V$ with $\|u_{0}\|+\|\theta_{0}\| \leq \rho,$
$\tilde{u}(t,\omega;t_{0},u_{0}) $ satisfies
\begin{eqnarray}\label{r-23}
\int_{-3}^{0}(|\nabla (|\tilde{u}|^{2})|^{2}+|\partial_{z}(|\tilde{u}|^{2})|^{2}+|\tilde{u}|\nabla \tilde {u}||^{2}
+|\tilde{u}|\partial_{z}\tilde{u}||^{2})dt \leq C_{2}(\omega).
\end{eqnarray}
By (\ref{qq-25}), (\ref{r-17}) and (\ref{r-23}), proceeding as (\ref{r-24}), for $\tilde{\mathbb{P}}-$a.e. $\omega\in \tilde{\Omega} $, there exists random variable $C_{3}(\omega)$ such that for arbitrary $\rho >0$ there exists $t(\omega)\leq -2$ such that for all $t_{0}\leq t(\omega)$ and $(u_{0}, \theta_{0})\in V$ with $\|u_{0}\|+\|\theta_{0}\| \leq \rho,$
$\bar{u}(t,\omega;t_{0},\bar{u}_{0}) $ satisfies
\begin{eqnarray}\label{r-26}
\sup\limits_{t\in [-2,0]} |\nabla\bar{u}(t,\omega; t_{0}, \bar{u}_{0} )|^{2}\leq C_{3}(\omega).
\end{eqnarray}
In view of (\ref{qq-31}) and (\ref{r-26}), following the steps in (\ref{r-24}), for $\tilde{\mathbb{P}}-$a.e. $\omega\in \tilde{\Omega} $, there exists a random variable $r_{5}(\omega)$,  such that for arbitrary $\rho >0$ there exists $t(\omega)\leq -1$ such that for all $t_{0}\leq t(\omega)$ and $(u_{0}, \theta_{0})\in V$ with $\|u_{0}\|+\|\theta_{0}\| \leq \rho,$ satisfies
\begin{eqnarray}\label{r-27}
 \sup\limits_{t\in [-1,0]}|\partial_{z}{u}(t,\omega; t_{0}, u_{0} )|^{2}+ \int_{-1}^{0} |\nabla{u}_{z}(s,\omega; t_{0}, u_{0} )|^{2}ds\leq r_{5}(\omega).
\end{eqnarray}
Regarding (\ref{r-29}), (\ref{r-9}), (\ref{r-16}) and (\ref{r-27}), we repeat the procedures of deriving (\ref{r-24}) and (\ref{r-22}). For $\tilde{\mathbb{P}}-$a.e. $\omega\in \tilde{\Omega} $, there exists a random variable $r_{6}(\omega)$ such that for arbitrary $\rho >0$ there exists $t(\omega)\leq -1$ such that for all $t_{0}\leq t(\omega)$ and $(u_{0}, \theta_{0})\in V$ with $\|u_{0}\|+\|\theta_{0}\| \leq \rho,$
$u(t,\omega;t_{0},u_{0}) $ satisfies
\begin{eqnarray}\label{r-30}
\sup\limits_{t\in[-1,0]} |\nabla {u}(t,\omega; t_{0}, u_{0} )|^{2} +\int_{-1}^{0} |\Delta{u}(s,\omega; t_{0}, u_{0} )|^{2}ds  \leq r_{6}(\omega).
\end{eqnarray}
By (\ref{qq-40}), (\ref{r-9}), (\ref{r-27}) and (\ref{r-30}), proceeding as above, for $\tilde{\mathbb{P}}-$a.e. $\omega\in \tilde{\Omega} $, there exists a random variable $r_{7}(\omega)$ such that for arbitrary $\rho >0$ there exists $t(\omega)\leq -1$ such that for all $t_{0}\leq t(\omega)$ and $(u_{0}, \theta_{0})\in V$ with $\|u_{0}\|+\|\theta_{0}\| \leq \rho,$
$\theta(t,\omega;t_{0},\theta_{0}) $ satisfies
\begin{eqnarray*}
\sup\limits_{t\in [-1,0]} \| \theta(t,\omega; t_{0}, \theta_{0} )\|^{2}\leq r_{7}(\omega).
\end{eqnarray*}
\par
Now we are ready to prove the desired compact result.
Let $r(\omega)= r_{5}(\omega)+r_{6}(\omega)+r_{7}(\omega)$, then  $B(-1,r(\omega))$, the ball of center $0\in V$ and radius $r(\omega),$ is an absorbing set
 at time $-1$ for $(S(t,s;\omega))_{t\geq s,\omega\in \tilde{\Omega}}$. According to Theorem \ref{thm-1}, in order to prove the existence of the random attractor in the space $V$, we need to to construct a compact absorbing set at time $0$ in $V$. Let $\mathcal{B}$ be a bounded subset of $V$, set
\begin{eqnarray*}
 \mathcal{C}_{T}:=\Big{\{} \Big{(}A_{1}^{\frac{1}{2}}\upsilon,   A_{2}^{\frac{1}{2}}T\Big{)}\Big{|}(\upsilon(-1),T(-1))\in \mathcal{B},
 (\upsilon(t),T(t))=S(t,-1;\omega)(\upsilon(-1),T(-1)),t\in[-1,0]\Big{\}}.
\end{eqnarray*}
We claim that $\mathcal{C}_{T}$ is compact in $L^{2}([-1,0]; {H} ).$  Indeed,
the space $V_{1}\times V_{2}  \subset H_{1}\times H_{2} $ is compact as $V_{i}\subset H_{i}$ is compact. Let $(\upsilon(-1), T(-1))\in \mathcal{B} ,$ by the argument of step 2 in the proof of Theorem \ref{thm-1}, we have
\[
(A_{1}^{\frac{1}{2}}u, A_{2}^{\frac{1}{2}}\theta)\in L^{2}([-1,0];V_{1}\times V_{2}),\ \
(\partial_{t}A_{1}^{\frac{1}{2}}u,  \partial_{t}A_{2}^{\frac{1}{2}}\theta )\in L^{2}([-1,0];V_{1}'\times V_{2}').
\]
Therefore, we deduce the result from Lemma \ref{lem-1} with
\[
B_{0}=V_{1}\times V_{2},\ \ B=H_{1}\times H_{2},\ \ B_{1}=V_{1}'\times V_{2}'.
\]
Now, we aim to show that for any fixed $t\in (-1,0],\omega\in \tilde{\Omega}, S(t,-1;\omega) $ is a compact operator in $V $.  Taking any bounded sequences
$\{(\nu_{0,n}, \tau_{0,n} ) \}_{n\in \mathbb{N}}$ in $ \mathcal{B}$, for any fixed $t\in (-1,0]$, $\omega\in \tilde{\Omega},$ we devote to extracting a convergent subsequence from
$\{S(t,-1;\omega)( \nu_{0,n}, \tau_{0,n}) \}$.
Since $\{(A_{1}^{\frac{1}{2}}\upsilon, A_{2}^{\frac{1}{2}}T  ) \}\subset \mathcal{C}_{T} ,$ by Lemma \ref{lem-1}, there is a function $(\nu_{*},\theta_{*})\in  L^{2}([-1,0];V)$
and a subsequence of $\{S(t,-1;\omega)(\nu_{0,n},\tau_{0,n} ) \}_{n\in \mathbb{N}}$ still denoted by  $\{S(t,-1;\omega)(\nu_{0,n},\tau_{0,n} ) \}_{n\in \mathbb{N}}$ such that
\begin{eqnarray}\label{r-31}
\lim\limits_{n\rightarrow \infty}\int_{-1}^{0}\|S(t,-1;\omega)(\nu_{0,n},\tau_{0,n} )-(\nu_{*}(t),\theta_{*}(t))  \|^{2}dt=0.
\end{eqnarray}
By the measure theory, we know that the convergence in mean square implies almost sure convergence. Therefore, it follows from (\ref{r-31}) that there exists a subsequence $\{S(t,-1;\omega)(\nu_{0,n},\tau_{0,n} ) \}_{n\in \mathbb{N}}$  such that
\begin{eqnarray}\label{r-32}
 \lim\limits_{n\rightarrow \infty}\|S(t,-1;\omega)(\nu_{0,n},\tau_{0,n} )-(\nu_{*}(t),\theta_{*}(t))  \|=0,\ \ a.e.\ t\in (-1,0].
\end{eqnarray}
Fix any $t\in (-1,0].$ By (\ref{r-32}), we can select a $t_{0}\in (-1,t)$ such that
\[
 \lim\limits_{n\rightarrow \infty}\|S(t_{0},-1,\omega)(\nu_{0,n},\tau_{0,n} )-(\nu_{*}(t_{0}),\theta_{*}(t_{0}))  \|=0.
\]
Then by the continuity of  $S(t-t_{0},t_{0};\omega) $ in $ V$ with respect to the initial value, we have
\begin{eqnarray*}
S(t,-1;\omega)(\nu_{0,n},\tau_{0,n})&=&S(t-t_{0},t_{0};\omega)S(t_{0},-1;\omega)(\nu_{0,n},\tau_{0,n})\\
&&\rightarrow S(t-t_{0},t_{0};\omega)(\nu_{*}(t_{0}),\theta_{*}(t_{0})),\ \ \ \mathrm{in}\ V.
\end{eqnarray*}
Hence, for any $t\in(-1,0]$, we can always find a convergent subsequence of $ \{S(t,-1;\omega)(\nu_{0,n},\tau_{0,n} ) \}_{n\in \mathbb{N}}$  in $V,$ which implies that for any fixed $t\in (-1,0],\omega\in \tilde{\Omega}, S(t,-1;\omega) $ is a compact operator in $V.$  Set
\[
\mathcal{B}(0,\omega)= \overline{S(0,-1;\omega)B(-1,r(\omega))},
\]
then, $\mathcal{B}(0,\omega)$ is a closed set of $ S(0,-1;\omega)B(-1,r(\omega))$ in $V.$ Using the above argument, we know  that $\mathcal{B}(0,\omega)$ is a random compact set in $V.$  Precisely,  $\mathcal{B}(0,\omega)$ is a compact absorbing set in $V$ at time $0.$ Indeed, for $(\nu_{0,n},\tau_{0,n} )\in \mathcal{B}, $ there exists $s(\mathcal{B})\in \mathbb{R}_{-}$ such that for any $s\leq s(\mathcal{B}),$ we have
\begin{eqnarray*}
S(0,s;\omega)(\nu_{0,n},\tau_{0,n})=S(0,-1;\omega)S(-1,s;\omega)(\nu_{0,n},\tau_{0,n})\subset S(0,-1;\omega) B(-1,r(\omega))\subset \mathcal{B}(0,\omega) .
\end{eqnarray*}
Therefore, we conclude the result from Theorem \ref{thm-1}.
\hspace{\fill}$\square$
\end{flushleft}

\section{Existence of invariant measure }
Up to now, we are ready to prove the existence of invariant measure of the system (\ref{e-1})-(\ref{e-5}).

Let $U_{0}=(v_{0},T_{0})\in  V$, $U(t,\omega; U_{0}):=(v(t,\omega;t_0,v_{0}), T(t,\omega;t_0,T_{0})) $ is the solution to (\ref{qq-1})-(\ref{qq-6}) with the initial value $U_0$. Following the standard argument, we can show that $U(t,\omega; U_{0}), t\in [t_0, \mathcal{T}]$ is Markov process in the sense that
for every bounded, $\mathcal{B}(V)$-measurable $F:V\rightarrow \mathbb{R},$ and all $s,t\in [t_{0}, \mathcal{T}]$, $t_{0}\leq s\leq t\leq  \mathcal{T}$,
\begin{eqnarray*}
\mathbb{E}(F(U(t,\omega; U_{0}))|\mathcal{F}_{s} )(\omega)=\mathbb{E}(F(U(t,s, U(s))) )\ \ \mathrm{for}\ \tilde{\mathbb{P}}-a.e.\ \omega\in \Omega,
\end{eqnarray*}
where $\mathcal{F}_{s}=\mathcal{F}_{t_{0},s}$ (see (\ref{e-13})), $ U(t,s, U(s))$ is the solution to (\ref{e-1})-(\ref{e-5}) at time $t$ with initial data $U(s).$
\par
For $B\in \mathcal{B}(V)$, define
\begin{eqnarray*}
\tilde{\mathbb{P}}_{t}(U_{0}, B)=\tilde{\mathbb{P}}((U(t,\omega; U_{0})\in B ).
\end{eqnarray*}
\par
For any probability  measure $\nu$ defined on $\mathcal{B}(V),$ denote the distribution at time $t$ of the solution to (\ref{e-1})-(\ref{e-5}) with initial distribution $\nu$ by
\[
(\nu \tilde{\mathbb{P}}_{t})(\cdot)= \int_{V}\tilde{\mathbb{P}}_{t}(x,\cdot)\nu(dx).
\]
\par
For $t\geq t_0$ and any continuous and bounded function $f\in C_{b}(V;\mathbb{R})$, we have
\begin{eqnarray*}
\tilde{\mathbb{P}}_{t}f(U_{0})=\mathbb{E}[f(U(t,\omega; U_{0})]=\int_{V}f(x)\tilde{\mathbb{P}}_{t}(U_{0}, dx).
\end{eqnarray*}
\begin{definition}
Let $\rho$ be a probability measure on $\mathcal{B}(V)$.  $\rho$ is called an invariant measure for $\tilde{\mathbb{P}}_{t}$, if
\begin{eqnarray*}
\int_{V}f(x)\rho(dx)=\int_{V}\tilde{\mathbb{P}}_{t}f(x)\rho(dx)
\end{eqnarray*}
for all $f\in C_{b}(V;\mathbb{R})$ and $t\geq 0.$
\end{definition}
\par
Let $\mu_{\cdot}$ be a transition probability from $\tilde{\Omega}$ to $V$, i.e., $\mu_{\cdot}$ is a Borel probability measure on $V$ and $\omega\rightarrow \mu_{\cdot}(B)$ is measurable for every Borel set $B\subset V.$ Denote by $\mathcal{P}_{\tilde{\Omega}}(V) $ the set of transition probabilities with $\mu_{\cdot}$ and $\nu_{\cdot}$ identified if $\tilde{\mathbb{P}}\{ \omega:\mu_{\omega}\neq \nu_{\omega}\}=0. $
\par
In view of Proposition $4.5$ in \cite{CF}, the existence of random attractor obtained in Theorem \ref{thm-3} implies the existence of invariant Markov measure $\mu_{\cdot}\in \mathcal{P}_{\tilde{\Omega}}(V)$ for $S$ such that $ \mu_{\omega}(\mathcal{A}(\omega))=1$, $ \tilde{\mathbb{P}}$-a.e.. Therefore, referring to $\cite{C}$,  there exists an invariant measure for the markov semigroup $\tilde{\mathbb{P}}_{t}$ and it is given by
$$\rho(B)=\int_{\tilde{\Omega}}\mu_{\omega}(B)\tilde{\mathbb{P}}(d\omega),$$
where $B\subseteq V$ is a Borel set. If the invariant measure $\rho$ for $\tilde {\mathbb{P}}$ is unique, the invariant Markov measure $\mu_{\cdot}$ for $S$ is unique and given by
$$\mu_{\omega}=\lim\limits_{t\rightarrow \infty}S(0,-t, \omega)\rho. $$
\par
Based on the above, we arrive at
\begin{theorem}
The Markov semigroup $ (\tilde{\mathbb{P}}_{t})_{t\geq 0}$ induced by the solution $(U(t,\omega; U_{0}))_{t\geq 0}$  to (\ref{e-1})-(\ref{e-5}) has an invariant measure $\rho$ with $ \rho(\mathcal{A}(\omega))=1\ \tilde{\mathbb{P}}$-a.e..
\end{theorem}



\def\refname{ Bibliography}


\begin{thebibliography}{2}

\bibitem {ARA}
R. A. Adams, Sobolev Spaces, Academic Press, New York, 1975.

\bibitem {C}
H. Crauel, Markov measures for random dynamical systems, Stochastics Stochastics
Rep.  3(1991), 153--173.

\bibitem {CDF}
H. Crauel, A. Debussche, F. Flandoli, Random attractors, J. Dynam. Differential Equations 9 (1997), 307--341.


\bibitem{CF}
H. Crauel,  F. Flandoli, Attractors  for  random dynamical  systems, Probab. Theory Relat. Fields. 100(1994), 365--393.


\bibitem {CIN}
C. Cao, S. Ibrahim, K. Nakanishi, E.S. Titi, Finite-time blowup for the inviscid primitive
equations of oceanic and atmospheric dynamics, Comm. Math. Phys. 337(2015), 473--482.




\bibitem {CLT1}
C. Cao, J. Li, E.S. Titi, Global well-posedness of strong solutions to the 3D primitive equations
with horizontal eddy diffusivity, J. Differential Equations 257 (2014), 4108--4132.


\bibitem {CLT2}
C. Cao, J. Li, E.S. Titi, Local and global well-posedness of strong solutions to the 3D primitive
equations with vertical eddy diffusivity, Arch. Ration. Mech. Anal. 214 (2014), 35--76.

\bibitem {CLT3}
C. Cao, J. Li, E.S. Titi, Global well-posedness of the three-dimensional primitive equations
with only horizontal viscosity and diffusion, Communications on Pure and Applied Mathematics  Vol. LXIX(2016), 1492--1531.


\bibitem {CT1}
C. Cao, E.S. Titi, Global well-posedness of the three-dimensional viscous primitive equations of large scale ocean and atmosphere dynamics, Ann. of Math. 166(2007), 245--267.


\bibitem {CT2}
C. Cao, E.S. Titi, Global well-posedness of the 3D primitive equations with partial vertical
turbulence mixing heat diffusion, Comm. Math. Phys. 310 (2012), 537--568.




\bibitem {DZ}
G. Da Prato,  and J. Zabczyk,   Stochastic equations in  infinite dimensions, Cambridge
University Press, Cambridge, 1992.

\bibitem{D-G-T-Z}  A. Debussche, N. Glatt-Holtz, R. Temam, M. Ziane: \emph{Global existence and regularity for the 3D stochastic primitive equations of the ocean and atmosphere with
multiplicative white noise}.  Nonlinearity, 25, 2093-2118 (2012).

  \bibitem{DZZ1}  Z. Dong, J. Zhai, R. Zhang, Large deviation principles for 3D stochastic primitive equations,  J. Differential Equations 263(2017), 3110-3146.


\bibitem{DZZ2}
Z. Dong, J. Zhai, R. Zhang, Exponential mixing for 3D stochastic primitive equations of the large scale ocean. preprint. Available at arXiv: 1506.08514.








\bibitem {G}
A. E. Gill, Atmosphere-ocean dynamics, International Geophysics Series, Vol. 30, Academic Press,
San Diego, 1982.




\bibitem {GH}
B. Guo  and D. Huang, 3d stochastic primitive equations of the large-scale ocean: global well-
posedness and attractors, Commun. Math. Phys. 286(2009), 697--723.






\bibitem {GKVZ}
N. Glatt-Holtz, I. Kukavica, V. Vicol, M. Ziane, Existence and regularity of invariant measures for
the three dimensional stochastic primitive equations, J. Math. Phys. 55(2014), 051504.


\bibitem {GMR}
F. Guill$\acute{\mathrm{e}}$n-Gonz$\acute{\mathrm{a}}$lez,  N. Masmoudi, M.A. Rodr$\acute{\mathrm{i}}$guez-Bellido, Anisotropic estimates and strong
solutions for the primitive equations, Diff. Int. Equ. 14(2001), 1381--1408.


\bibitem {H1}
G. J. Haltiner, Numerical weather prediction, J.W. Wiley \& Sons,  New York, 1971.

\bibitem {H2}
G. J. Haltiner and R. T. Williams, Numerical prediction and dynamic meteorology, John Wiley \&
Sons, New York, 1980.

\bibitem {HTZ}
 C. Hu, R. Temam, M. Ziane, The primimitive equations of the large scale ocean under the small depth
hypothesis, Disc. and Cont. Dyn. Sys. 9(2003), 97--131.


\bibitem {J}
N. Ju, The global attractor for the solutions to the 3D viscous primitive equations, Discret. Contin. Dyn.
Syst. 17(2007), 159--179.


\bibitem {KZ}
I. Kukavica  and M. Ziane, On the regularity of the primitive equations of the ocean, Nonlinearity
20(2007), 2739--2753.





\bibitem {LM}
J. Lions and B. Magenes, “Nonhomogeneous Boundary Value Problems and Applications,”
Springer-Verlag, New York, 1972.

\bibitem {LTW1}
J.L. Lions,  R. Temam and S. Wang, New formulations of the primitive equations of atmosphere and applications, Nonlinearity 5(1992), 237--288.

\bibitem {LTW2}
J.L. Lions,  R. Temam and S. Wang, On the equations of the large scale ocean, Nonlinearity 5(1992), 1007--1053.

\bibitem {LTW3}
J.L. Lions,  R. Temam and S. Wang, Models of the coupled atmosphere and ocean$(CAO I)$, Computational
Mechanics Advance 1(1993), 1--54.

\bibitem {LTW4}
J.L. Lions,  R. Temam and S. Wang, Mathematical theory for the coupled atmosphere-ocean models
$(CAO III)$, J. Math. Pures Appl. 74(1995), 105--163





\bibitem {T2}
R. Temam, “ Navier-Stokes equations. Theory and Numerical Analysis,” reprint of 3rd edition,
AMS 2001.




\bibitem{Z}
G. Zhou, Random attractor of the 3D viscous primitive equations driven by fractional noises, 	arXiv:1604.05376.





\end{thebibliography}
\end{document}